\documentclass[a4paper,12pt,leqno]{article}

\newcommand{\kmcomment}[1]{}
\usepackage{amsmath,amsfonts,amssymb}
\usepackage{comment}
\usepackage{mathtools}
\newcommand{\rank}{\operatorname{rank}}

\newcommand{\ovfrakh}{\overline{\frakh}}
\newcommand{\bfrakg}{\bar{\frakg}}

\newcommand{\ds}{\ensuremath{\displaystyle }}
\newcommand{\pdel}{\partial}

\newcommand{\myHom}[1]{\textrm{H}_{#1}}
\newcommand{\myHomW}[2]{\textrm{H}_{#1,#2}}
\newcommand{\mN}{\ensuremath{\mathbb{N}}} %Use in Math-mode
\newcommand{\mR}{\ensuremath{\mathbb{R}}} %Use in Math-mode
\newcommand{\mZ}{\ensuremath{\mathbb{Z}}} %Use in Math-mode
\newcommand{\xb}[1]{x_{#1}}
\newcommand{\yb}[1]{y_{#1}}
\newcommand{\zb}[1]{z_{#1}}
\newcommand{\wb}[1]{w_{#1}}
\newcommand{\ndeg}[1]{\operatorname{ndeg}(#1)}

\newcommand{\frakg}{\mathfrak{g}}

\newcommand{\frakh}{\mathfrak{h}}
\newcommand{\frakhN}[1]{\mathfrak{h}_{#1}}

\newcommand{\tbdl}[1]{\mathrm{T}(#1)}

\newcommand{\cbdl}[1]{\mathrm{T}^{*}(#1)}

\newcommand{\Sbt}[2]{[#1,#2]}                
\newcommand{\SbtD}[2]{[\thickspace[#1,#2]\thickspace]}  % Double               
\newcommand{\SbtDE}[3]{{[[#2,#3]]}_{#1}}  % Double               
\newcommand{\SbtS}[2]{[#1,#2]_{\text{\footnotesize S}}} %Use in Math-ode 
                
\newcommand{\SbtE}[3]{\ensuremath{[#2,#3]_{#1}}} % Enhanced                
\newcommand{\parity}[1]{(-1)^{#1}}

\newcommand{\Lb}[1]{\mathit{L}_{#1}}  % Lie bibun (differentiation)
\newcommand{\inn}[1]{i_{#1}} 
\newcommand{\TCkt}[2]{[\![#1,#2]\!]}

\newcommand{\SbtX}[3]{\ensuremath{[#2,#3]_{#1}}}                
                
\newcommand{\Pkt}[2]{\{#1,#2\}}% 

\newcommand{\BktD}[3]{\{#2,#3\}_{#1}}

\newcommand{\cgaiseki}[1]{\Lambda^{#1} \cbdl{M}}

\newcommand{\CS}[2]{ \text{C}_{#2,#1}} % chain space 

\newcommand{\CSp}[1]{\text{C}_{#1}} % chain space 

\renewcommand{\SbtD}{\TCkt}

\newcommand{\SbtSZ}[2]{[#1 , #2]_{sz}} % like SbtSS

\newcommand{\mw}{\mywedge}
\newcommand{\we}{\wedge}

\newcommand{\wfg}[1]{ \frakg_{<#1>}} % double frakg 

\renewcommand{\dim}{\textrm{dim}}

\newcommand{\myCSW}[2]{ \text{C}_{#1,#2}} 
\newcommand{\myOSW}[2]{ \overline{\text{C}}_{#1,#2}} 
\newcommand{\mywedge}{\bigtriangleup} 

\newcommand{\myCS}[1]{ \text{C}_{#1}} % chain space

\renewcommand{\[}{$$} \renewcommand{\]}{$$}

\renewcommand{\SbtD}{\TCkt}

\usepackage{wrapfig}

\usepackage{fancyvrb,shortvrb}
 \RequirePackage{fancyvrb}
  \DefineVerbatimEnvironment%
      {newverb}{Verbatim}
      {xleftmargin=10mm}
 
 \DefineVerbatimEnvironment%
     {newverb*}{Verbatim}
     {xleftmargin=10mm,showspaces=true,showtabs=true}

 \usepackage{listings}
  \RequirePackage{listings}
  \lstnewenvironment{kmverb}
     { \lstset{xleftmargin=10mm, breaklines=true}}
     { }

  \lstnewenvironment{kmverbb}[1]
     { \lstset{xleftmargin=#1, showspaces=false, breaklines=true}}
     { }
  \lstnewenvironment{kmverbN10}
     { \lstset{xleftmargin=-10mm, breaklines=true}}
     { }
  \lstnewenvironment{kmverbZero}
     { \lstset{xleftmargin=0mm, breaklines=true}}
     { }
 
  \lstnewenvironment{kmverb*}
     { \lstset{xleftmargin=10mm, showspaces=true, breaklines=true}}
     { }
 
\usepackage[a4paper]{geometry}
\usepackage{amsthm}
\theoremstyle{definition} %plain, definiton, remark の三種類を指定
\newtheorem{theorem}{Theorem}
\newtheorem*{theorem*}{Theorem}
\newtheorem{definition}[theorem]{Definition}
\newtheorem*{definition*}{Definition}
\newtheorem{example}{Example}[section]
\newtheorem{prop}{Proposition}[section]
\newtheorem{remark}{Remark}
\newtheorem{lemma}{Lemma}

\title{Superalgebra structure on differential forms of manifold }
\author{Kentaro Mikami \and Tadayoshi Mizutani}

\usepackage{setspace}
\usepackage[a4paper]{geometry}
\usepackage[svgnames]{xcolor}

\AtBeginDocument{
  \abovedisplayskip     =0.5\abovedisplayskip
  \abovedisplayshortskip=0.5\abovedisplayshortskip
  \belowdisplayskip     =0.5\belowdisplayskip
  \belowdisplayshortskip=0.5\belowdisplayshortskip}

\numberwithin{equation}{section}

\allowdisplaybreaks
\AtBeginDocument{\parindent=0pt \setlength{\headheight}{16pt}}
\begin{document}
\maketitle
\tableofcontents
\section{Introduction}

Typical examples of Lie superalgebra are the direct sum of exterior power of
tangent bundle of a differential manifold with the Schouten bracket or the
matrix algebra whose elements are  divided into 4 parts, and the artificial
bracket operation.  The Schouten bracket in the first example is
useful to describe a 2-vector field  \(\pi\) to be a  
Poisson tensor by \(\Sbt{\pi}{\pi} = 0\).  

In this article, we introduce a notion of Lie superalgebra structure on 
the direct sum of exterior power of
cotangent bundle of a differential manifold with a natural bracket.  
Also we introduce an extension of the Lie superalgebra structure by the
1-vector fields.  

\subsection{Quick review of \(\mZ\)-graded Lie superalgebras}
First we recall the definition of Lie superalgebra or pre Lie
superalgebra we call sometimes.   
\begin{definition}[\(\mZ\)-graded Lie superalgebra]
Suppose a real vector space 
$\frakg $ is graded by \(\ds \mZ\) as 
\(\ds \frakg = \sum_{j\in \mZ} \frakg_{j} \)
and has a \(\mR\)-bilinear operation \(\Sbt{\cdot}{\cdot}\)  satisfying 
\begin{align}
& \Sbt{ \frakg_{i}}{ \frakg_{j}} \subset \frakg_{i+j} \label{cond:1} \\
& \Sbt{X}{Y} + (-1) ^{ x y} 
 \Sbt{Y}{X} = 0 \quad \text{ where }  X\in \frakg_{x} \text{ and } 
 Y\in \frakg_{y}  \\
& 
(-1)^{x z} \Sbt{ \Sbt{X}{Y}}{Z}  
+(-1)^{y x} \Sbt{ \Sbt{Y}{Z}}{X}  
+(-1)^{z y} \Sbt{ \Sbt{Z}{X}}{Y}  = 0\;.  
%\quad \text{(Jacobi identity).} 
\label{super:Jacobi}
\end{align}
Then we call \(\frakg\) a \textit{\(\mZ\)-graded} or \textit{(pre) Lie superalgebra}.    
\kmcomment{
A Lie  \textit{superalgebra} $\frakg $ is graded by \(\ds \mZ_{2}\) as \(\ds
\frakg = \frakg_{[0]} \oplus \frakg_{[1]} \) and the condition
\eqref{cond:1} is regarded as  \(\ds \Sbt{ \frakg_{[1]}}{ \frakg_{[1]}}
\subset \frakg_{[0]}\) in modulo 2 sense.   
}
\end{definition}

\begin{remark} 
Super Jacobi identity \eqref{super:Jacobi} above is equivalent to the one of
the following.  
\begin{align}
\Sbt{X}{\Sbt{Y}{Z}} &=  
\Sbt{\Sbt{X}{Y}}{Z}  + (-1)^{x y} \Sbt{Y}{ \Sbt{X}{Z}} 
\\  
\Sbt{ \Sbt{X}{Y}}{Z} &=  
\Sbt{X}{ \Sbt{Y}{Z} } + (-1)^{y z} \Sbt{ \Sbt{X}{Z}}{Y} 
\end{align} 
\kmcomment{
Suppose \(\ds \frakg = \sum_{j\in \mZ} \frakg_{j} \) is a pre Lie
superalgebra.  Let \(\ds\frakg_{[0]} = \sum_{i \text{ is even}}
\frakg_{i}\) and \(\ds\frakg_{[1]} = \sum_{i \text{ is odd}}
\frakg_{i}\). Then \(\ds \frakg = \frakg_{[0]} \oplus \frakg_{[1]} \)
holds and this is a Lie superalgebra.  
}
\end{remark}

\kmcomment{
\medskip

\begin{example} \label{exam:super:elem} 
Take an $n$-dimensional vector space $V$ and split it as \(\ds V =
V_{0} \oplus V_{1}\).  Define \(\ds \frakg_{[i]} = \{ A\in
\mathfrak{gl}(V) \mid A (V_{j}) \subset V_{i+j} \} \), where 
\(\ds \mathfrak{gl}(V)\) is the space of endomorphisms of \(V\). 
For each \(\ds A \in \frakg_{[i]} \) and \(\ds B \in \frakg_{[j]} \), 
define \(\ds\Sbt{A}{B} = A B - (-1)^{i j} B A \). 
Then  \(\ds \mathfrak{gl}(V)
=\frakg_{[0]} \oplus \frakg_{[1]} \) with this bracket is a Lie
superalgebra. 

More concretely, we take $n=2$ and \(\ds \dim V_{[i]}=1\) for \(i=0,1\).  
Then \(\ds \frakg_{[0]} = \begin{bmatrix} \ast & 0 \\ 0 & \ast
\end{bmatrix} \) and 
\(\ds \frakg_{[1]} = \begin{bmatrix}  0 & \ast  \\ \ast & 0
\end{bmatrix} \).    
If we define 
\(\ds \frakg_{0} =  \begin{bmatrix} a & 0 \\ 0 & -a \end{bmatrix}\), 
\(\ds \frakg_{1} = \begin{bmatrix}  0 & \ast  \\ \ast & 0
\end{bmatrix}\) and 
\(\ds \frakg_{2} =  \begin{bmatrix} a & 0 \\ 0 & a \end{bmatrix}\). Then 
\(\ds \mathfrak{gl}(2) =\frakg_{0}\oplus\frakg_{1}\oplus\frakg_{2}\) is
a pre Lie superalgebra.  
\end{example}

\medskip
}

\kmcomment{
We will introduce the notion of \textit{double weight} in pre Lie superalgebras
(cf.\ Definition \ref{defn:w:weight}) and our results in this note are
the calculation of the Euler number and Betti numbers of homology groups of
double-weighted pre  Lie superalgebras of special type.  

(cf.\ 
Lemma \ref{thm:Euler:wZero},  
Theorem \ref{thm:triv:general}  
For general $n$, the Euler number of chain complex 
\(\ds \{ \myCSW{\bullet}{w,h}\} \) 
is 0 for 
each $w$ and each $h$.  

Theorem \ref{thm:module:gen}). 
}%endOFkmcomment

\begin{definition}[Sub superalgebra]
For a general \(\mZ\)-graded Lie superalgebra  
\(\ds \frakg = \sum_{i} \frakg_{i}\),   
 if \(\frakh_{i}\) is a subspace of \(\frakg_{i} \) for each \(i\), 
 and satisfy \( \Sbt{\frakh_{i}} {\frakh_{j}} \subset \frakh_{i+j}\),
 then \(\ds \sum_{i} \frakh_{i}\) is a \(\mZ\)-graded Lie
 superalgebra, we call a Lie 
 \textcolor{red}{sub}-superalgebra of \(\frakg\).  
 \end{definition}
\begin{remark}
Let \(\ds \frakg = \sum_{i} \frakg_{i}\) be a \(\mZ\)-graded Lie
 superalgebra. 
From definitions, we see some obvious properties.  
\begin{itemize}
\item[(1)]
 Assume \(\frakg_{\ell} \ne (0)\). Then
\(\ds \bfrakg^{\ell} = \sum_{i>0} \frakg_{i\ell}\) or  
\(\ds \sum_{i \geqq 0} \frakg_{i\ell} = \frakg_{0} + \bfrakg^{\ell}
\)   
 are \(\mZ\)-graded Lie \textcolor{red}{sub}-superalgebras of
 \(\frakg\).  If \(\ell\) is even, they are Lie algebras.  
\item[(2)]
 Fix \( \ell \geqq 0\).
  \(\ds \sum_{i \geqq \ell}\frakg_{ i}\) is a  
 \(\mZ\)-graded Lie \textcolor{red}{sub}-superalgebra of \(\frakg\). 

\item[($2'$)]
 Fix \( \ell \leqq 0\).
 \(\ds \sum_{i \leqq \ell}\frakg_{ i}\) is a  
 \(\mZ\)-graded Lie \textcolor{red}{sub}-superalgebra of \(\frakg\). 
\end{itemize}
\end{remark}

%NHK

\section{Super brackets on differential forms}
We introduce a new example of \(\mZ\)-graded Lie superalgebra with a super
bracket associated from the exterior differentiation. 
We have a typical \(\mZ\)-graded Lie superalgebra of the direct sum of  \(j\)
-multivectors \(\Lambda^{j} \tbdl{M}\) whose super grade (or weight) is
\(j-1\) for \( j = 0, \ldots, \dim M\), where  \(M\) is a manifold. 
\[
\begin{array}{|c|*{8}{c}|}
\hline
 \Lambda^{\bullet} \tbdl{M} &&&&& \text{1-vec} & \text{ 2-vec} &\text{ 3-vec}
 & \cdots \\ \hline
grading & \cdots &-3 & -2 & -1& 0 & 1 & 2 & \cdots \\ \hline
\Lambda^{\bullet}\cbdl{M} &\cdots & 2\text{-form} & 1\text{-form} &
0\text{-form}&   &   &   & \cdots \\ 
\hline
\end{array}
\]
So, the super grade of differential \(i\)-forms \(\Lambda^{i} \cbdl{M}\) of
``\(\mZ\)-graded Lie superalgebra'' of the direct sum of differential 
\(i\)-forms \(\Lambda^{i} \cbdl{M}\) are expected to be \( -(i+1)\), and we
define the super grade (weight) of differential forms as follows:

\begin{definition} \label{defn:grading}
The grade (in super sense) or  weight of 
 \( \alpha \in \Lambda^{a} \cbdl{M}\) is defined by \( -(1+a)\). Sometimes we
 abbreviate the grade  \(-(1+a) \) of \( \alpha \in \Lambda^{a} \cbdl{M}\) 
 as  \(a' \).  
\end{definition}

\medskip

\begin{definition}
We define a bi-linear map 
\begin{equation}
\SbtD{\alpha} {\beta} = \parity{a} d( \alpha \wedge \beta )
= \alpha \wedge d \beta - \parity{a' b'} \beta \wedge d \alpha
\label{defn:bkt:forms}
\end{equation}
for \( \alpha \in \Lambda^{a} \cbdl{M} \)
and \( \beta \in \Lambda^{b} \cbdl{M} \). 
\end{definition}

\begin{theorem}
The direct sum \(\ds \Lambda^{n}\cbdl{M} \oplus \cdots  \oplus \Lambda^{1}
\cbdl{M} \oplus \Lambda^{0}\cbdl{M} \) becomes  a \(\mZ\)-graded Lie
superalgebra with the bracket \eqref{defn:bkt:forms} and the grading by
Definition \ref{defn:grading}. 

\end{theorem}
\textbf{Proof:}
Suppose 
\( \alpha \in \Lambda^{a} \cbdl{M} = \frakg_{[a']}\) and 
\( \beta  \in \Lambda^{b} \cbdl{M} = \frakg_{[b']}\), where \( a'=-(1+a)\)
and \( b' = - (1+b)\).  
Then \( \SbtD{\alpha}{\beta} \in \Lambda^{1+a+b} \cbdl{M} = \frakg_{[a'+b']}
\) because  \( -1 - (1+a+b) = a' + b'\). 
\begin{align*}
\SbtD{\beta}{\alpha} &= \beta\wedge d \alpha - \parity{b' a'}
\alpha\wedge\beta 
= \parity{1+a'b'} ( \alpha \wedge 
 d\beta - \parity{a' b'} \beta\wedge d \alpha 
)
=  \parity{1+a'b'} \SbtD{\alpha}{\beta}
\\
\parity{a' c'} \SbtD{\alpha}{\SbtD{\beta}{\gamma}} &= 
\parity{a' c'}(   \alpha \wedge d (\SbtD{\beta}{\gamma} ) 
- \parity{ a'(b'+c') } \SbtD{\beta}{\gamma} \wedge d\alpha )
\\&= 
\parity{a' c'}(  0 
- \parity{ a'(b'+c') }
( \beta \wedge d\gamma - \parity{b'c'}\gamma\wedge d\beta  ) 
\wedge d\alpha )
\\&= 
- \parity{a' c'}   
 \parity{ a'(b'+c') } \beta \wedge d\gamma \wedge d\alpha 
+ 
\parity{a' c'}   
\parity{ a'(b'+c') } \parity{b'c'}\gamma\wedge d\beta \wedge d \alpha
\\&= 
- \parity{a' b'}   
 \beta \wedge d\gamma \wedge d\alpha 
+ 
 \parity{b'c'} \gamma\wedge  d \alpha \wedge d\beta
 \\\shortintertext{and so}
 \mathfrak{S}_{\alpha,\beta,\gamma} & 
\parity{a' c'} \SbtD{\alpha}{\SbtD{\beta}{\gamma}} = 0 \;.
\end{align*}
\qed

\begin{remark}
We remember that in the superalgebra \(\sum_{i\geqq 0} \Lambda^{i} \tbdl{M}
\), \( \SbtS{f}{g} =0\) because of \(\frakg_{-2} = (0)\). 
But the definition above says  if \(\alpha,\beta\) are 0-forms, namely
functions \(f,g\), then \( \SbtD{f}{g} = d( f g )\), which is a 1-form.  
\end{remark}

\begin{remark}
It is a great surprise to the authors if this theorem was previously unknown.

Concerning to choice of bracket, we do not claim that \eqref{defn:bkt:forms}
is unique. In fact, a common constant multiple of \eqref{defn:bkt:forms} is
also super bracket.  \end{remark}

\begin{remark}
A constant 1 satisfies \(\SbtD{1}{1} = 0 \) and \( \phi : A \mapsto
\SbtD{1}{A}\) is a derivation of degree \(1\), and satisfies \(\phi \circ \phi = 0\)
like the Poisson cohomology story.  In fact, this is the coboundary operation  of de Rham cohomology theory. 
\end{remark}

Some values of super bracket are 
\begin{align*}
\SbtD{1}{\alpha} &= \parity{a} \SbtD{\alpha}{1} = d \alpha \;,
\\
\SbtD{1}{\alpha\wedge\beta} & = d(\alpha\wedge \beta)= \parity{a}
\SbtD{\alpha}{\beta}\;,  
\\
\SbtD{1}{\alpha\wedge\beta} & = d(\alpha\wedge \beta)
= (d\alpha)\wedge \beta + \parity{a} \alpha\wedge d \beta 
= \SbtD{1}{\alpha} \wedge \beta + 
\parity{a} \alpha \wedge \SbtD{1}{\beta}\;.  
\\\shortintertext{More generally}
\SbtD{\gamma}{\alpha\wedge\beta} & = 
 \SbtD{\gamma}{\alpha} \wedge \beta + 
\parity{c' (1+a')} \alpha \wedge \SbtD{\gamma}{\beta}
+ \parity{c'}(d\gamma) \wedge  \alpha \wedge \beta
\;, \\\shortintertext{because}  
\SbtD{\gamma}{\alpha\wedge\beta} 
& = \parity{c}d ( \alpha \wedge \beta \wedge \gamma ) 
 = \parity{c} (d \alpha) \wedge \beta \wedge \gamma  
+  \alpha \wedge( d \beta) \wedge \gamma  
+  \parity{a}  \alpha \wedge \beta \wedge (d  \gamma ) 
\\& 
= \SbtD{\gamma}{\alpha} \wedge \beta + \parity{a+ac}
   \alpha \wedge \gamma  \wedge d \beta 
\\& = \SbtD{\gamma}{\alpha} \wedge \beta + \parity{a+ac}
   \alpha \wedge \left( \SbtD{\gamma}{\beta} - \parity{c}(d \gamma) \wedge
   \beta\right) 
\\& = \SbtD{\gamma}{\alpha} \wedge \beta + \parity{a+ac}
   \alpha \wedge  \SbtD{\gamma}{\beta} + \parity{c+1}(d \gamma) \wedge\alpha
   \wedge \beta \;.  
\end{align*}

\begin{remark} We may give a new grading  by the grade of \(i\)-form by
\(1+i\) if we only deal with the direct sum of forms, but we expect to
handle differential forms and multi-vector fields, too.  
\end{remark}

\section{Homology groups of our Lie superalgebra}
 
\subsection{Quick review of homology groups of general Lie superalgebra}
%\section{Preliminaries, Notations and Basic Facts} \label{sec:prelim} 
In a usual Lie algebra homology theory, $m$-th chain space is the
exterior product \(\ds \Lambda^{m}\frakg\) of \(\frakg\) and the
boundary operator essentially comes from the operation \(\ds X \wedge
Y \mapsto \Sbt{X}{Y}\). 

Likewise, in the case of \(\mZ\)-graded Lie superalgebras, 
"exterior algebra" is defined as the quotient of the tensor
algebra \(\ds \otimes^{m}  \frakg \) of \(\frakg \) by the two-sided
ideal generated by  \begin{align} & X \otimes Y + (-1)^{x y} Y \otimes X
\quad\text{where }\quad X \in \frakg_{x}, Y \in \frakg_{y} \;,
\end{align} and we denote the equivalence class of \(\ds X \otimes Y \)
by \(\ds X \mywedge Y\).  
\kmcomment{ \(\ds \bigtriangleup \triangle
\mywedge  \nabla \bigtriangledown \) } 
Since \( \ds X_{\text{odd}}
\mywedge Y_{\text{odd}}  = Y_{\text{odd}} \mywedge X_{\text{odd}} \) and
\( \ds X_{\text{even}} \mywedge Y_{\text{any}}  = - Y_{\text{any}}
\mywedge X_{\text{even}} \) hold, \( \mywedge ^{m} \frakg_{k} \) has a
symmetric property for odd $k$ and has a skew-symmetric property 
for even $k$ with respect to \(\mywedge\).   

\kmcomment{
\medskip
\begin{definition}
Assume that the pre Lie superalgebra \(\ds \frakg\) acts on a module
$V$ as follows: For each homogeneous \(\xi \in \frakg\) 
there corresponds an element 
\(\ds \xi_{V} \in \text{End}(V) \) and satisfies  
\(\ds \Sbt{\xi}{\eta}_{V} = \xi_{V}  \circ \eta_{V}  
- (-1)^{|\xi| |\eta| } \eta_{V} \circ \xi_{V} \) where 
\(\ds \xi \in \frakg_{|\xi|}\) and 
\(\ds \eta \in \frakg_{|\eta|}\). Then we call $V$ a \(\frakg\)-module.   

\kmcomment{
In short, 
we have a pre Lie superalgebra homomorphism  
\(\ds \frakg \ni \xi \mapsto \xi_{V} \in \text{End}(V)\).  
}%endOFkmcomment
We often write \(\ds \xi\cdot v\) for \(\ds \xi_{V}(v) \).   
\end{definition}
\begin{example} A pre Lie superalgebra \(\ds \frakg\) is itself
a \(\frakg\)-module by own bracket
\(\ds X \cdot Z = \Sbt{X}{Z}\).  
Let $X,Y\in\frakg$ be homogeneous as \(\ds X\in \frakg_{x}, 
 Y\in \frakg_{y}\).   Then,  
\(\ds (X\circ  Y - (-1)^{x y} Y\circ X ) \cdot Z = 
\Sbt{X}{Y}\cdot Z \) holds and this is just Jacobi identity.  
\end{example}
}

% \subsection{Homology groups weighted by the first grading} 
\begin{definition}
For a \(\mZ\)-graded Lie superalgebra 
\(\ds \frakg = \sum_{i} \frakg_{i} \), $m$-th chain space is defined by 
\(\CSp{m} = \underbrace{\frakg  \mw  \cdots  \mw  \frakg}_{m}\)
and the boundary operator is given by 
\begin{align} 
\pdel (A_{1}  \mw  A_{2}  \mw  
\cdots  \mw  A_{m}) 
=& \sum_{i<j} \parity{ i-1 + a_{i}\sum_{i<s<j} a_{s} }
A_{1}  \mw  \cdots 
\widehat{ A_{i} }\cdots  \mw  
\SbtD{A_{i}}{A_{j}}  \mw  
 \cdots  \mw  A_{m}
 \label{bdary:orig}
\end{align}
where \( A_{i} \in \frakg_{a_{i}}\). We see that \( \pdel \circ \pdel = 0\)
holds, and have the homology groups by 
\[
\myHom{m} (\frakg ) =  \ker(\pdel_{} : \myCS{m} \rightarrow
\myCS{m-1})/ \pdel ( \myCS{m+1} )\;.  
\]
 
%and \( \pdel \CSp{m} \subset \CSp{m-1}\) 
\kmcomment{
\begin{equation}
\CS{w}{m} = \sum_{i_{1}\leqq \cdots\leqq i_{m}}\{\frakg_{i_{1}}  \mw 
\cdots  \mw  \frakg_{i_{m}} \mid i_{1}+\cdots+i_{m}= w\}\;. 
\end{equation}
}

\end{definition}

\begin{definition}[Weight]
For each non-zero element \(u\) in \(\ds \frakg_{i_{1}} \mywedge  \cdots \mywedge
\frakg_{i_{m}} \), we define  \(u\) has the weight 
\(\ds i_{1}+\dots + i_{m} \). 
We 
define the subspace of 
\(\myCS{m} \) by 
\[\ds \myCSW{m}{w} = 
 \sum_{
\substack{
i_{1}\leq \ldots \leq i_{m}\\
\sum_{s=1}^{m}i_{s} = w
}
} 
\frakg_{i_{1}} \mywedge  \cdots \mywedge
\frakg_{i_{m}} \;,\] which is the space of  common weight $w$.  
\end{definition}
The following is known well. 
\begin{prop}
The weight $w$ is preserved by \(\ds \pdel\), i.e., we have 
\(\ds \pdel( \myCSW{m}{w} ) \subset \myCSW{m-1}{w} \). Thus, 
for a fixed $w$, we have
$w$-weighted homology groups
\[\ds 
\myHomW{m}{w} (\frakg ) =  \ker(\pdel : \myCSW{m}{w} \rightarrow
\myCSW{m-1}{w})/ \pdel ( \myCSW{m+1}{w} )\;.  
\] 
\end{prop}

Since the boundary operator is defined by 
\eqref{bdary:orig}, we see 
\begin{align} 
\pdel (A_{1}  \mw  
\cdots  \mw  A_{p} \wedge 
B_{1}  \mw 
\cdots  \mw  B_{q} 
) &= 
\pdel (A_{1}  \mw  \cdots  \mw  A_{p})  \wedge B_{1}  \mw  \cdots
 \mw  B_{q}  %\notag
\\&\quad 
+ \parity{p} A_{1}  \mw  \cdots  \mw  A_{p}  
\wedge  \pdel (B_{1}  \mw  \cdots  \mw  B_{q} )  
\notag
\\&\quad   
+ \SbtSZ{ A_{1}  \mw  \cdots  \mw  A_{p} }
{ B_{1}  \mw  \cdots  \mw  B_{q} }
\notag
\\\shortintertext{where}
 \SbtSZ{ A_{1}  \mw  \cdots  \mw  A_{p} }
{ B_{1}  \mw  \cdots  \mw  B_{q} }
&= \sum
 \parity{ i-1 + a_{i}(\sum_{i<s}a_{s} + \sum_{s<j} b_{s})} 
\cdots \widehat{A_{i}}
 \cdots  \mw  A_{p}  
\wedge  B_{1}  \mw  \cdots \SbtD{A_{i}}{B_{j}} \cdots 
\label{sbt:sz}
\\
&= \sum
 \parity{ p+j + (\sum_{i<s}a_{s} + \sum_{t<j} b_{t}) b_{j}} 
\cdots \SbtD{A_{i}}{B_{j}}
 \cdots  \mw  A_{p}  
\wedge  B_{1}  \mw  \cdots \widehat{B_{j}} \cdots 
\\&= \sum \parity{i+ a_{i}\sum _{s>i} a_{s} 
+ j +b_{j}\sum _{s<j} b_{s}} \cdots \widehat{ A_{i} } \cdot\cdot A_{p}  \mw  
\SbtD{A_{i}}{B_{j}} \mw  B_{1} \cdot\cdot  \widehat{ B_{j} } \cdots \; .
\end{align}
This new bracket 
\(\SbtSZ{\cdot}{\cdot}\), which is similar to the Schouten bracket in some
case,   satisfies 
\begin{align*}
\SbtSZ{ A_{1} \cdots A_{p} }
{ B_{1} \cdots B_{q} \mw  C_{1} \cdots C_{r} }
&= 
\SbtSZ{ A_{1} \cdots A_{p} }
{ B_{1} \cdots B_{q}} \mw  C_{1} \cdots C_{r} 
\\&\quad 
+\parity{ 1+q + (\sum_{i}a_i)( \sum_{j} b_j ) }  B_{1} \cdots B_{q} \mw 
\SbtSZ{ A_{1} \cdots A_{p} }
{  C_{1} \cdots C_{r} }
\\\shortintertext{and also}
\SbtSZ{ A_{1} \cdots A_{p} 
\mw  B_{1} \cdots B_{q} }{ C_{1} \cdots C_{r} }
&= \parity{p} A_{1} \cdots A_{p} \mw 
\SbtSZ{ B_{1} \cdots B_{q}} { C_{1} \cdots C_{r} }  
\\&\quad 
+\parity{ 1+ (\sum_{j}b_j)( \sum_{k} c_k ) }  
\SbtSZ{ A_{1} \cdots A_{p} }
{  C_{1} \cdots C_{r} }
\mw 
B_{1} \cdots B_{q} \; . 
\end{align*}
Assume that \( A_{1} = \cdots = A_{p} = A \in \frakg_{\text{odd}}\) in \eqref{sbt:sz}.
Then we have
\begin{align} \SbtSZ{  \mw ^{p}A }{B_{1}  \mw  \cdots  \mw  B_{q}}
& = p \parity{p}  \mw ^{p-1} A  \mw \sum_{j}\parity{j+ b_{j} \sum_{s<j}b_{s}} \sum_{j=1}^{q} 
\SbtD{ A }{ B_{j} }  \mw B_{1}  \cdots
\\& 
 =
 %\textcolor{red}{-} 
 p   \mw ^{p-1}(- A) \mw \SbtSZ{A}{{B_{1} \mw \cdots \mw  B_{q}} } \; . 
\\\shortintertext{Furthermore, assume  every \(B_j\) has degree \textcolor{red}
{even}, then }
\SbtSZ{ A }{B_{1}  \mw  \cdots  \mw  B_{q}} &= \sum_{j}
 B_{1} \mw \cdots \mw
\SbtSZ{ A }{B_{j}} \mw  \cdots  \mw  B_{q} \; .
\\\shortintertext{Contrary, assume  every \(C_j\) has degree \textcolor{red}
{odd}, then }
\SbtSZ{ A }{C_{1}  \mw  \cdots  \mw  C_{r}} &= \sum_{j}\parity{j+1}
 C_{1} \mw \cdots \mw
\SbtSZ{ A }{C_{j}} \mw  \cdots  \mw  C_{r}\\& = 
\sum_{j}  
\SbtSZ{ A }{C_{j}} \mw  C_{1} \mw \cdots \mw \widehat{C_{j}} \mw \cdots  \mw  C_{r}  
\\\shortintertext{and }
 \SbtSZ{ A }{B_{1}  \mw  \cdots  \mw  B_{q}
\mw
C_{1}  \mw  \cdots  \mw  C_{r}}
 &= 
\SbtSZ{ A }{B_{1}  \mw  \cdots  \mw  B_{q}} \mw 
C_{1}  \mw  \cdots  \mw  C_{r} \\&\quad 
+ B_{1}  \mw  \cdots  \mw  B_{q}  \mw 
\SbtSZ{ A }{C_{1}  \mw  \cdots  \mw  C_{r}} 
\notag \; .
\end{align}
In particular, 
the boundary operator is written by left action as follows: 
\begin{align} 
\pdel (A_{0}  \mw  A_{1}  \mw  
\cdots  \mw  A_{m}) 
=& 
- A_{0}  \mw  \pdel (A_{1}  \mw  \cdots  \mw  A_{m})
+ A_{0} \cdot (A_{1}  \mw  \cdots  \mw  A_{m})
\label{pdel:recurs:left}
\\\shortintertext{where} 
 A_{0} \cdot (A_{1}  \mw  \cdots  \mw  A_{m}) 
=& \sum_{i=1}^{m} (-1)^{ a_{0}\sum_{1\leqq s<i}a_{s}} 
A_{1}  \mw  \cdots  \mw  \SbtD{A_{0}}{A_{i}}  \mw  \cdots 
 \mw  A_{m} 
 \\
=&  - \SbtSZ{ A_{0} } { A_{1}  \mw  \cdots  \mw  A_{m} }
 \label{left:action}
\end{align}
for each homogeneous elements \( A_{i}\in\frakg_{a_{i}}\).

In lower degree, the boundary operator is given
as below: 
\begin{align} 
\pdel ( A  \mw  B ) &= \SbtD{A}{B} \\
\pdel ( A  \mw  B  \mw  C ) &= - {A}  \mw  {
\SbtD{B}{C} } + 
\SbtD{A}{B}  \mw  C + (-1)^{a b} B  \mw  
\SbtD{A}{C} 
\end{align}
for each homogeneous elements 
\( A\in\frakg_{a}\), \( B\in\frakg_{b}\), \( C\in\frakg_{c}\).

\subsection{Some works of (co)homology groups of \(\mZ\)-graded Lie
superalgebra with the Schouten bracket}
We only refer to \cite{Mik:Miz:super2} and  \cite{Mik:Miz:super3}.  

\kmcomment{
\section{A new example of \(\mZ\)-graded Lie superalgebra with a super
bracket associated from the exterior differentiation}
}

\section{Homology groups of superalgebra of left invariant Lie groups}
Assume a Lie group \(G\) acts on \(M\). Then we are able to discuss 
 \(G\)-invariant theory of \(\mZ\)-graded Lie superalgebra of tangent bundle or
cotangent bundle. In particular case of \(M=G\), we already have studied  
of \(\mZ\)-graded Lie superalgebra of Lie algebra, where  
let \( \xi_{i}\) be a basis of Lie algebra with bracket relation 
\begin{align}
\Sbt{\xi_{i}}{\xi_{j}} &= \sum_{k} c^{k}_{ij} \xi_{k} \\
\shortintertext{with 
\(  c^{k}_{ij}\) are constants of structure. 
Let \( \zb{i}\) be the dual
of \( \xi_{j}\). Since} 
(d \zb{i})(\xi_{j}, \xi_{k}) &= 
\xi_{j} \langle \zb{i}, \xi_{k} \rangle 
- \xi_{k} \langle \zb{i}, \xi_{j} \rangle 
- \langle \zb{i} , \Sbt{\xi_{j}}{\xi_{k}} \rangle 
= - c^{i}_{jk}\notag \\
d \zb{i} &= \sum_{j,k}
(d \zb{i})(\xi_{j}, \xi_{k})  \zb{j} \otimes \zb{k}
= - \sum_{j,k} c^{i}_{jk} \zb{j} \otimes \zb{k}\notag \\
&= -\frac{1}{2} \sum_{j,k} c^{i}_{jk} \zb{j} \otimes \zb{k}
 -\frac{1}{2} \sum_{j,k} c^{i}_{kj} \zb{k} \otimes \zb{j}
= -\frac{1}{2} \sum_{j,k} c^{i}_{jk}
(\zb{j} \otimes \zb{k} -\zb{k} \otimes \zb{j})
\notag
 \\\shortintertext{we have }
 d \zb{i} &= 
 - \textcolor{red}{ \frac{1}{2}}\sum_{j,k} c^{i}_{jk} \zb{j} \wedge \zb{k} =
 - \sum_{j<k} c^{i}_{jk} \zb{j} \wedge \zb{k} \;,
\label{ext:deriv:d}
\end{align}

\subsection{When \(M=G\)  and \(\dim G=2\)}
The \(w\)-weighted \(m\)-th chain space \(\CS{w}{m}\) is given by 
\[ \sum_{i_{1},i_{2},i_{3}} \frakg_{-1}^{i_{1}}  \mw  
 \frakg_{-2}^{i_{2}}  \mw  
 \frakg_{-3}^{i_{3}} 
\text{ where }  i_{1} + i_{2} + i_{3} = m\;  \text{ and } 
 i_{1} + 2 i_{2} + 3 i_{3} = - w\;. \]  
We see that \( m \leqq -w\) and \( 3 m \geqq -w\), namely 
\( Q(m) \leqq 0  \), where \( Q(m) =  (-w -m ) (-w -3m)  \).    
\(m\) runs from 
\( \text{ceil}( (-w)/3 ) \) to \(-w\). 
Since the dimension of \(\frakg_{-2}\) is 2 and \(  \mw ^{i_{2}}
\frakg_{-2}\) is a skew-symmetric power, and so it vanishes if \( i_{2}>2\).
Thus, solving linear equations of \( i_{1}, i_{3} \) for cases of \(i_{2}=0,
1,2\), we see that 

\begin{align}
\CS{w}{m} & = 
\begin{cases}
 \mw ^{-w-3K} \frakg_{-1} % \mw  
    \mw ^{0}\frakg_{-2} % \mw  
    \mw ^{K}\frakg_{-3} 
+ 
    \mw ^{-w-3K-1} \frakg_{-1} % \mw  
    \mw ^{2}\frakg_{-2} % \mw  
    \mw ^{K-1}\frakg_{-3} 
   & \text{ if } -w -m = 2 K
   \\
    \mw ^{-w-3L-2} \frakg_{-1} % \mw  
    \mw ^{1}\frakg_{-2} % \mw  
    \mw ^{L}\frakg_{-3} 
   & \text{ if } -w -m = 2 L + 1 \;,
   \end{cases}
   \notag
\\\shortintertext{in short} 
&= \begin{cases}
 \mw ^{-w-3K} 1 % \mw  
    \mw ^{K} V 
+ 
    \mw ^{-w-3K-1} 1  \mw  
   \zb{1}  \mw  \zb{2} % \mw  
    \mw ^{K-1} V 
   & \text{ if } -w -m = 2 K
   \\
    \mw ^{-w-3L-2} 1 % \mw  
    \mw ^{1}\zb{i} % \mw  
    \mw ^{L}V 
   & \text{ if } -w -m = 2 L + 1 \;. 
\end{cases}
\label{cs:basis}
\end{align}

For lower weight \(w\), the chain complexes are as below. 
\[
\begin{array}{|c||*{7}{c|}}
\hline
w &\CS{w}{1}&\CS{w}{2}&\CS{w}{3}&\CS{w}{4}&\CS{w}{5}&\CS{w}{6}\\ 
\hline
\hline
-1 & \frakg_{-1} & 0  & 0  & 0  & 0  & 0 
\\ \hline
-2 & \frakg_{-2} &  \mw ^{2}\frakg_{-1}   & 0  & 0  & 0  & 0 \\
\hline
-3 & \frakg_{-3} & \frakg_{-2} \mw \frakg_{-1}& \mw ^{3}\frakg_{-1}& 0 & 0  & 0 
\\ \hline
-4 & 0 & \frakg_{-3} \mw \frakg_{-1} +  \mw ^{2}\frakg_{-2}
&
\frakg_{-2}  \mw   \mw ^{2}\frakg_{-1} & \mw ^{4}\frakg_{-1} & 0 & 0 
\\ \hline
-5 & 0 &  \frakg_{-3} \mw \frakg_{-2} 
& \frakg_{-3} \mw  \mw ^{2}\frakg_{-1} +  \mw ^{2}\frakg_{-2}
 \mw  \frakg_{-1} &
\frakg_{-2}  \mw   \mw ^{3}\frakg_{-1} & \mw ^{5}\frakg_{-1} & 0 
\\ \hline
\end{array}
\]
\kmcomment{  
When \(w=-2\), the second chain space \( \CS{-2}{2}\) is 1-dimensional and
spanned by \( 1    \mw  1\) and by \( \pdel ( 1  \mw  1) = 0\).  
When \(w=-3\), we need \( \pdel (\zb{i}  \mw  1) = - d \zb{i}\).  

For \(w=-4\), we used 
 \( \pdel ( \zb{i}  \mw  1  \mw  1)   
 = - \zb{i}  \mw  \SbtD{1}{1}
 + \SbtD{ \zb{i}}{1}  \mw  1
 +\parity{(-1)(-1)} 1  \mw   \SbtD{ \zb{i}}{1}
= 2 ( 
  \SbtD{ \zb{i}}{1}  \mw  1
) = -2 ( 1  \mw  d \zb{i} ) \), 
 \( \pdel ( (\zb{i} \wedge \zb{j})  \mw  1)   
 = \SbtD{\zb{i} \wedge  \zb{j}}{1} = 0\),  
and 
\( \pdel (\zb{i}  \mw  \zb{j}) = 
\SbtD{\zb{i}} { \zb{j}} = 0\).  
  
For \(w=-5\), we used 
 \( \pdel(( \zb{1} \wedge  \zb{2})   \mw  \zb{j} )   
 = \SbtD{ \zb{1} \wedge  \zb{2}}{\zb{j}}  = 0\),   

\( \pdel( \zb{1}  \mw   \zb{2}   \mw  1 )   = -   \zb{1}
 \mw  \SbtD{ \zb{2}}{ 1 } + \SbtD{\zb{1}}{\zb{2}}  \mw  1
+ \zb{2}  \mw  \SbtD{ \zb{1} }{ 1 }  = \zb{2}  \mw  ( - d
\zb{1} ) = 2 \zb{2}  \mw  (  \zb{1} \wedge \zb{2})  \),  

\( \pdel(( \zb{1} \wedge  \zb{2})   \mw  1  \mw  1 ) = 0\),   
\( \pdel ( \zb{i}  \mw  1  \mw  1  \mw  1)   = 0  + \zb{i}
\cdot (  1  \mw  1  \mw  1 ) = 3 (  1  \mw  1  \mw 
\SbtD{\zb{i}}{1} ) = -3  (  1  \mw  1  \mw   d \zb{i} ) \).  
}%endOFkmcomment
\medskip

%nhk

The dimension of 
 each chain space is directly given by 
\begin{align*}
\dim 
\CS{w}{m} 
&= \begin{cases}
1 & \text{ if } Q(m) = 0 
\\
2 & \text{ if } -w -m = 2 K > 0  \text{ and } Q(m) < 0 
\\
2 & \text{ if } -w -m = 2 L + 1 \text{ and }  Q(m) < 0 
   \\
0 & \text{ otherwise } 
\end{cases}
% \\\shortintertext{and } &
= \begin{cases}
1 & \text{ if } Q(m) = 0  
\\
2 & \text{ if }  Q(m) < 0  
\\
0 & \text{ otherwise}   
   \;. 
\end{cases}
\end{align*} 
Non-abelian Lie algebra is unique and the
structure equation is \( \Sbt{\xi_{1}}{\xi_{2}} = \xi_{1}\) and so
\( d \zb{1} = - \zb{1} \wedge \zb{2}\),  
\( d \zb{2} = 0\). 
We put \(\zb{1}\wedge \zb{2}\) by \(V\). Then \( d V = 0\). 

The multiplications of super bracket are shown in the table below. 
\[
\begin{array} {|c||*{4}{c|}} \hline
 & 1\in \frakg_{-1} & \zb{i} \in \frakg_{-2} 
 & V = \zb{1}\wedge \zb{2} \in \frakg_{-3} 
 \\\hline\hline
1\in \frakg_{-1} & 0 & d \zb{i} 
 &  0 
 \\\hline
\zb{p} \in \frakg_{-2}  & - d \zb{p} & 
0 
& 0  \\\hline
V = \zb{1}\wedge \zb{2} \in \frakg_{-3} &  0  
& 0 & 0  \\\hline 
\end{array}
\]

We see that 
\begin{align}
& 
\pdel (  \mw ^{a} 1  \mw ^{c} V ) = 0 \text{
because of } \SbtD{1}{1} = 0,  \SbtD{1}{V} = \SbtD{V}{1} = 0\text{  and
}  \SbtD{V}{V} = 0\;. 
\label{dImage:cases:1}
\\  
& 
\pdel (  \mw ^{a} 1  \mw  \zb{1}  \mw \zb{2}  \mw ^{c}
V )  = a \zb{2}  \mw ^{a-1} 1  \mw ^{c+1} V 
\label{dImage:cases:2}
\\
& 
\pdel (  \mw ^{a} 1  \mw  \zb{i}   \mw ^{c}
V )  = \parity{a}a  \mw ^{a-1} 1  \mw (\zb{i}\cdot 1)  \mw ^{c} V
= \begin{cases} 
\parity{a} a  \mw ^{a-1} 1   \mw ^{c+1} V & \text{ if } i=1 \\
0 & \text{ if } i=2
\end{cases}
\label{dImage:cases:3}
\end{align} 
because of 
\begin{align*}
\pdel( \mw ^{a} 1  \mw  \zb{1} \mw \zb{2} \mw ^{c}V)
&= 
\pdel(\zb{1} \mw \zb{2} \mw ^{a} 1  \mw ^{c}V)
= - \zb{1}  \mw  \pdel (
\zb{2} \mw ^{a} 1  \mw ^{c}V)
+ \zb{1}\cdot ( \zb{2} \mw ^{a} 1  \mw ^{c}V)
\\ \pdel (
\zb{2} \mw ^{a} 1  \mw ^{c}V)
&= - \zb{2}  \mw  \pdel (
 \mw ^{a} 1  \mw ^{c}V)
+ \zb{2}\cdot ( \mw ^{a} 1  \mw ^{c}V) = 0 
\\
\zb{1}\cdot ( \zb{2} \mw ^{a} 1  \mw ^{c}V)
 &= 
( \zb{1}\cdot  \zb{2})  \mw ^{a} 1  \mw ^{c}V
 + a \zb{2} \mw (\zb{1}\cdot 1) \mw ^{a-1}1   \mw ^{c}V
 + c \zb{2} \mw ^{a} 1  \mw  (\zb{1}\cdot V)   \mw ^{c-1}V
\\
 &= 
0 
 + a \zb{2} \mw (  V ) \mw ^{a-1}1   \mw ^{c}V
 + 0 
 =  a \zb{2} \mw ^{a-1}1   \mw ^{c+1}V
 \; . 
\\
\pdel (  \mw ^{a} 1  \mw  \zb{i}   \mw ^{c} V ) 
&= \parity{a}
\pdel (  \zb{i}  \mw ^{a} 1    \mw ^{c} V ) 
 = \parity{a}a  \mw ^{a-1} 1  \mw (\zb{i}\cdot 1)  \mw ^{c} V
 \; . 
\end{align*}

From \eqref{dImage:cases:1} \(\sim\) \eqref{dImage:cases:3}, \eqref{cs:basis}
shows generators of the boundary image of each chain space as follows: 
\begin{align}
\pdel \CS{w}{m} 
&= \begin{cases}
0  & \text{ if } Q(m) = 0 
\\
  (-w-3K-1) \zb{2}  \mw ^{-w-3K-2} 1  \mw ^{K} V 
   & \text{ if } Q(m) < 0 \text{ and } -w -m = 2 K 
   \\
   \parity{-w-3L-2}  (-w-3L-2)  \mw ^{-w-3L-3} 1  \mw ^{L+1}V 
   & \text{ if } Q(m) < 0 \text{ and } -w -m = 2 L + 1 \\
   0 & \text{ otherwise }\;. 
\end{cases}
\label{dImage:cs:basis}
\\\shortintertext{and so}
\dim (\pdel \CS{w}{m} )  
&= \begin{cases}
0  & \text{ if } Q(m) = 0 \text{ or } m = \text{ceil}(-w/3) 
\\
 1 & \text{ if } Q(m) < 0 \text{ and }  m > \text{ceil}(-w/3) 
 \\
0 & \text{ otherwise } 
\end{cases}
= \begin{cases}
 1 & \text{ if } Q(m) < 0 \text{ and }  m > \text{ceil}(-w/3) 
 \\
0 & \text{ otherwise. } 
\end{cases}
\end{align}
In short, both ends of chain complex vanish by the boundary operator
\(\pdel\). Each chain space between both ends has rank 1. Thus, we have the
following tables of space dimensions and rank of the boundary operator and
the Betti numbers: When \( -w = 3 \Omega + \varepsilon\) with \(
\varepsilon=\pm 1\) and we put \( \Omega_{0} = \text{ceil}( -w/3 ) = 
\Omega + (\varepsilon + 1)/2 \). 

\begin{center}
\(\begin{array} { c | *{1}{c} }
-w=1 & 1  \\
\hline
\dim &  1\\
\rank & 0 \\
\hline
Betti &  1
\end{array}\)
\hfil
\(\begin{array} { c | *{2}{c} }
-w=2 & 1 & 2 \\
\hline
\dim & 2 & 1\\
\rank & 0 & 0 \\
\hline
Betti & 2 & 1
\end{array}\)
\hfil
%\( -w = 3 \Omega\)-case: 
\(
\begin{array}{c| *{7}{c}}
%\text{degree}
 -w = 3 \Omega 
& \Omega & \Omega + 1 &\Omega+2 &  \cdots & 3\Omega -2 &  3\Omega -1 &
3\Omega \\\hline 
SP\dim & 1 & 2 & 2 & \cdots & 2 & 2 & 1
\\
\rank & 1 & 1 & 1 & \cdots & 1 & 0 & 0 \\\hline
Betti & 0 & 0 & 0 & \cdots & 0 & 1 & 1 \\\hline
\end{array}
\)

\medskip

%\( -w = 3 \Omega+ 1\)-case: 
\(
\begin{array}{c| *{7}{c}}
%\text{degree}
 -w = 3 \Omega+ \varepsilon > 3
& \Omega_{0} & \Omega_{0} + 1 &\Omega_{0}+2 &  \cdots & 3\Omega+
\varepsilon-2  &
3\Omega + \varepsilon-1 &
3\Omega+\varepsilon \\\hline 
SP\dim & 2 & 2 & 2 & \cdots & 2 & 2 & 1
\\
\rank & 1 & 1 & 1 & \cdots & 1 & 0 & 0 \\\hline
Betti & 1 & 0 & 0 & \cdots & 0 & 1 & 1 \\\hline
\end{array}
\)

\end{center}

\kmcomment{
%\( -w = 3 \Omega+ 2\)-case: 

\(
\begin{array}{c| *{7}{c}}
%\text{degree}
-w = 3 \Omega+ 2  
& \Omega+1 & \Omega + 2 &\Omega+3 &  \cdots & 3\Omega  &  3\Omega+1  &
3\Omega+2 \\\hline 
SP\dim & 2 & 2 & 2 & \cdots & 2 & 2 & 1
\\
\rank & 1 & 1 & 1 & \cdots & 1 & 0 & 0 \\\hline
Betti & 1 & 0 & 0 & \cdots & 0 & 1 & 1 \\\hline
\end{array}
\)
}

\kmcomment{
\\ =
\dim (\pdel \CS{w}{m} )  
&= \begin{cases}
0  & \text{ if } Q(m) = 0 
\\
 1 & \text{ if } -w -m = 2 K > 0 \text{ and } -w-3K-2 \geqq 0 
\\
 0 & \text{ if } -w -m = 2 K > 0  \text{ and } -w-3K-2 < 0 
\\
1 & \text{ if } -w -m = 2 L + 1 \text{ and } -w-3L-3  \geqq 0 \\
0 & \text{ if } -w -m = 2 L + 1 \text{ and } -w-3L-3 < 0   
\end{cases}
\\& 
= \begin{cases}
1 &\text{ if } 0 < -w -m = 2 K \text{ and } 0<K\leqq \frac{1}{3}( -w-2) 
\\
1 & \text{ if } -w -m-1 = 2 L \text{ and } L \leqq  \frac{1}{3} ( -w-3) \\
0 & \text{ otherwise } 
\end{cases}
\\& 
= \begin{cases}
1 &\text{ if } Q(m) < 0 \text{ and } m \neq  \text{ceil}(w/3)   
\\
0 & \text{ otherwise } \;. 
\end{cases}
\end{align}

\begin{center}
\(\begin{array} { c | *{2}{c} }
-w=2 & 1 & 2 \\
\hline
\dim & 2 & 1\\
\rank & 0 & 0 \\
\hline
Betti & 2 & 1
\end{array}\)
\hfil 
\(\begin{array} { c | *{3}{c} }
-w=3 & 1 & 2 & 3 \\
\hline
\dim & 1 &  2 & 1\\
\rank & 1 & 0 & 0 \\
\hline
Betti & 0 & 1 & 1
\end{array}\)
\hfil 
\(\begin{array} { c | *{4}{c} }
-w=4   & 2 & 3 & 4 \\
\hline
\dim &   2 &  2 & 1\\
\rank &  1 & 0 & 0 \\
\hline
Betti  & 1 & 1 & 1
\end{array}\)
\hfil 
\(\begin{array} { c | *{5}{c} }
-w=5  & 2 & 3 & 4 & 5 \\
\hline
\dim  & 2&   2 &  2 & 1\\
\rank & 1 &   1 & 0 & 0 \\
\hline
Betti & 1 & 0 & 1 & 1
\end{array}\)
\hfil 
\(\begin{array} { c | *{5}{c} }
-w=6  & 2 & 3 & 4 & 5 & 6 \\
\hline
\dim & 1 & 2& 2 &  2 & 1\\
\rank & 1& 1 &   1 & 0 & 0 \\
\hline
Betti & 0 & 0 & 0 & 1 & 1
\end{array}\)
\end{center}

}%endOFkmcomment

\bigskip

\subsection{When \(M=G\)  and \(\dim G=3\)} 
\subsubsection{Common issues}

``General'' multiplication table of 
super bracket is shown as follows: 
\[
\begin{array} {|c||*{4}{c|}} \hline
 & 1\in \frakg_{-1} & \zb{i} \in \frakg_{-2} 
 & \zb{j}\wedge \zb{k} \in \frakg_{-3} 
 &  \zb{1}\wedge   \zb{2}\wedge   \zb{3} \in \frakg_{-4} 
 \\\hline\hline
1\in \frakg_{-1} & 0 & d \zb{i} 
 & d( \zb{j}\wedge \zb{k}) 
 &  0  
 \\\hline
\zb{p} \in \frakg_{-2}  & - d \zb{p} & 
- d(\zb{p}\wedge \zb{i}) 
& 0 & 0 \\\hline
\zb{q}\wedge \zb{r} \in \frakg_{-3} & d(\zb{p}\wedge \zb{q}) 
& 0 & 0 & 0 \\\hline 
\zb{1}\wedge   \zb{2}\wedge   \zb{3}
 \in \frakg_{-4} & 0 &0 & 0 & 0  
 \\\hline
\end{array}
\]

\kmcomment{
We use the notation \(\zb{i}\) instead of \(\zb{i}\) in
\eqref{ext:deriv:d} and 
}
We denote the 2-form \( \zb{i} \wedge \zb{i+1}\) by \(
\wb{i+2}\), where indices are reduced by modulo 3.
Put \( \zb{1} \wedge \zb{2} \wedge \zb{3}\) by \(V\). 

Each 
chain space \(\CS{w}{m}\) consists of spaces \( 
 \mw ^{i_{1}}\frakg_{-1}
 \mw ^{i_{2}}\frakg_{-2}
 \mw ^{i_{3}}\frakg_{-3}
 \mw ^{i_{4}}\frakg_{-4}
\) 
where \(i_{2}=0,1,2,3\) and \(i_{4}=0,1\), 
and 
\begin{align}
m &= i_{1}+ i_{2} + i_{3} + i_{4} \quad \text{
and }\quad - w =  i_{1}+ 2 i_{2} + 3 i_{3} + 4 i_{4} \;.   
\label{weight:cond}
\end{align}
From \(
  i_{1}+ i_{2} + i_{3} = m - i_{4}\) and \(
  i_{1}+ 2 i_{2} + 3 i_{3} = - w - 4 i_{4} \),    
we see \( m - i_{4} \leqq -w - 4 i_{4} \leqq 3( m - i_{4} )\), i.e., 
\begin{equation}
 \frac{- w-1}{3} \leqq m \leqq -w \;.  
 \end{equation}
Solving the linear equations \eqref{weight:cond}, we have 
\begin{alignat*}{3}
2 i_{1} &= 3 m + w - i_{2} + i_{4} \;, 
& 
2 i_{3} & = - w  - m - i_{2} - 3 i_{4} \;. 
\\\shortintertext{Depending on the cases \(-w-m = 2K\) or 
 \(-w-m = 2L+1\), the above equations become}
2 i_{1} &= - 2 w - 6 K - i_{2} + i_{4} \;, 
& 
2 i_{3} & = 2 K  - i_{2} - 3 i_{4} 
\\
2 i_{1} &= - 2 w - 6L - 3 -  i_{2} + i_{4} \;, 
& \quad 
2 i_{3} &= 2 L +1  - i_{2} - 3 i_{4} \;.
\end{alignat*}
When 
\( -w = 3\Omega + \varepsilon\), \(\varepsilon=0,\pm 1\), 
\begin{align*}
2 i_{1} & = - 2 w - 6 K - i_{2} + i_{4} 
=  2 (3\Omega + \varepsilon)  - 6 K - i_{2} + i_{4} 
\;, 
\  
2 i_{3} = 2 K  - i_{2} - 3 i_{4} \quad 
 \text{ if } i_{2}+ i_{4} \text{ is even, and }     
\\
2 i_{1} & = - 2 w - 6L - 3 -  i_{2} + i_{4} 
 =  2 ( 3\Omega + \varepsilon )  - 6L - 3 -  i_{2} + i_{4} \;, 
 \  
2 i_{3} = 2 L +1  - i_{2} - 3 i_{4} \quad 
\text{ if } i_{2}+ i_{4} \text{ is odd.}  
\end{align*}
Thus, 
%\begin{subequations}
\begin{align}
\CS{-w-2K}{w} &= \sum_{i_{2}+i_{4}= even} 
 \mw ^{ 3\Omega + \varepsilon -3K - i_{2}/2+i_{4}/2 } \frakg_{-1}
 \mw ^{ i_{2} } \frakg_{-2}
 \mw ^{ K - i_{2}/2- 3 i_{4}/2 } \frakg_{-3}
 \mw ^{ i_{4} } \frakg_{-4}
\label{G3:codim:even}
\\&= 
 \mw ^{ - w  -3K } \frakg_{-1}
 \mw ^{0} \frakg_{-2}
 \mw ^{ K } \frakg_{-3}
 \mw ^{ 0 } \frakg_{-4}
+
 \mw ^{ - w -3K - 1 } \frakg_{-1}
 \mw ^{ 2 } \frakg_{-2}
 \mw ^{ K  - 1 } \frakg_{-3}
 \mw ^{ 0 } \frakg_{-4}
\notag
\\&\quad 
+
 \mw ^{ - w -3K  } \frakg_{-1}
 \mw ^{ 1 } \frakg_{-2}
 \mw ^{ K - 2 } \frakg_{-3}
 \mw ^{ 1 } \frakg_{-4}
+
 \mw ^{ - w  -3K -1 } \frakg_{-1}
 \mw ^{ 3 } \frakg_{-2}
 \mw ^{ K - 3 } \frakg_{-3}
 \mw ^{ 1 } \frakg_{-4}
\notag
\\
\CS{-w -2L-1}{w} &= \sum_{i_{2}+i_{4}= odd} 
 \mw ^{ 3\Omega+\varepsilon -3L -3/2- i_{2}/2+i_{4}/2 } \frakg_{-1}
 \mw ^{ i_{2} } \frakg_{-2}
 \mw ^{ L + 1/2 - i_{2}/2- 3 i_{4}/2 } \frakg_{-3}
 \mw ^{ i_{4} } \frakg_{-4}
\label{G3:codim:odd}
\\&= 
 \mw ^{-w  -3L -2} \frakg_{-1}
 \mw ^{1} \frakg_{-2}
 \mw ^{ L } \frakg_{-3}
 \mw ^{ 0 } \frakg_{-4}
+
 \mw ^{-w  -3L - 3 } \frakg_{-1}
 \mw ^{ 3 } \frakg_{-2}
 \mw ^{ L  - 1 } \frakg_{-3}
 \mw ^{ 0 } \frakg_{-4}
\notag
\\&\quad 
+
 \mw ^{-w -3L -1 } \frakg_{-1}
 \mw ^{ 0 } \frakg_{-2}
 \mw ^{ L - 1 } \frakg_{-3}
 \mw ^{ 1 } \frakg_{-4}
+
 \mw ^{-w  -3L -2 } \frakg_{-1}
 \mw ^{ 2 } \frakg_{-2}
 \mw ^{ L - 2 } \frakg_{-3}
 \mw ^{ 1 } \frakg_{-4}
\; . \notag
\end{align}
%\end{subequations}

The expressions above imply the following dimension formula, where
\(\tbinom{p}{q}\) is \( \frac{ p! }{q! (p-q)!}\) only when \(p\geqq 0\), 
\( q\geqq 0\) and \( p\geqq q\). \(\tbinom{p}{q} \) is 0 if   
\( p<0\) or \(q<0 \) or \( p <  q\). 
\begin{subequations}
\begin{align} 
\dim \CS{w}{-w-2K} 
&= \tbinom{1-1+3\Omega+\varepsilon- 3 K}{1-1} \tbinom{3-1+ K}{3-1}
+\tbinom{1-1+3\Omega+\varepsilon-3 K-1}{1-1}\tbinom{3}{2}\tbinom{3-1+ K-1}{K-1}
\\&\quad 
+\tbinom{1-1+3\Omega+\varepsilon-3 K}{1-1}\tbinom{3}{1}\tbinom{3-1+ K-2}{K-2}
+\tbinom{1-1+3\Omega+\varepsilon-3 K-1}{1-1}\tbinom{3}{3}\tbinom{3-1+ K-3}{K-3}
\notag
\\ &= \tbinom{3\Omega+\varepsilon- 3 K}{0} \tbinom{K+2}{2}
+3 \tbinom{3\Omega+\varepsilon-3 K-1}{0}\tbinom{ K+1}{K-1}
\\&\quad 
+3 \tbinom{3\Omega+\varepsilon-3 K}{0}\tbinom{ K}{K-2}
+\tbinom{3\Omega+\varepsilon-3 K-1}{0}\tbinom{ K-1}{K-3}
\notag
\\ &= \tbinom{3\Omega+\varepsilon- 3 K}{0}(\tbinom{K+2}{2}
+3 \tbinom{ K}{K-2} ) 
%\\&\quad 
+ \tbinom{3\Omega+\varepsilon-3 K-1}{0}( 3  \tbinom{ K+1}{K-1}
+\tbinom{ K-1}{K-3} ) 
\\&= \begin{cases}
\tbinom{K+2}{2}
+ 3 \tbinom{ K+1}{K-1}
+ 3 \tbinom{ K}{K-2}  
+ \tbinom{ K-1}{K-3}  
& \text{ if } K \leqq \Omega + \frac{\varepsilon-1}{3} \\
\tbinom{K+2}{2} +3 \tbinom{ K}{K-2}  
& \text{ if } \varepsilon=0 \text{ and } K = \Omega 
\\
0 & \text{ otherwise }
\end{cases} 
\end{align}
\begin{align} 
\dim \CS{w}{-w-2L-1} 
&= \tbinom{1-1+3\Omega+\varepsilon- 3 L-2}{1-1}\tbinom{3}{1} \tbinom{3-1+ L}{3-1}
+\tbinom{1-1+3\Omega+\varepsilon-3 L-3}{1-1}\tbinom{3}{3}\tbinom{3-1+ L-1}{L-1}
\\&\quad 
+\tbinom{1-1+3\Omega+\varepsilon-3 L-1}{1-1}\tbinom{3-1+ L-1}{L-1}
+\tbinom{1-1+3\Omega+\varepsilon-3 L-2}{1-1}\tbinom{3}{2}\tbinom{3-1+ L-2}{L-2}
\notag
\\ &= \tbinom{3\Omega+\varepsilon- 3 L-2}{0} 3 \tbinom{L+2}{2}
+\tbinom{3\Omega+\varepsilon-3 L-3}{0}\tbinom{L+1}{L-1}
\\&\quad 
+\tbinom{3\Omega+\varepsilon-3 L-1}{0}\tbinom{ L+1}{L-1}
+\tbinom{3\Omega+\varepsilon-3 L-2}{0}3 \tbinom{L}{L-2}
\notag
\\ &=3  \tbinom{3\Omega+\varepsilon- 3 L-2}{0} 
\left( \tbinom{L+2}{2} + \tbinom{L}{L-2} \right)
%\\&\quad 
+\left(\tbinom{3\Omega+\varepsilon-3 L-3}{0}
+\tbinom{3\Omega+\varepsilon-3 L-1}{0}\right)
    \tbinom{ L+1}{L-1}
\\& = \begin{cases}
3 ( \tbinom{L+2}{2} + \tbinom{L}{L-2} )
+ 2 \tbinom{ L+1}{L-1}
& \text{ if } L \leqq \Omega + \frac{\varepsilon}{3}-1 \\
\tbinom{ L+1}{L-1}
+ 3 ( \tbinom{L+2}{2} + \tbinom{L}{L-2} )
& \text{ if } \varepsilon= -1 \text{ and } L = \Omega -1  
\\
  \tbinom{ L+1}{L-1}
& \text{ if } \varepsilon= 1 \text{ and } L = \Omega   
\\
0 & \text{ otherwise.}
\end{cases} 
\end{align}
\end{subequations}

For lower weight \(w\), the chain complexes are as below. 
\[
\begin{array}{|c||*{7}{c|}}
\hline
w &\CS{w}{1}&\CS{w}{2}&\CS{w}{3}&\CS{w}{4}&\CS{w}{5}&\CS{w}{6}\\ 
\hline
\hline
-1 & \frakg_{-1} & 0  & 0  & 0  & 0  & 0 
\\ \hline
-2 & \frakg_{-2} &  \mw ^{2}\frakg_{-1}   & 0  & 0  & 0  & 0 \\
\hline
-3 & \frakg_{-3} & \frakg_{-1} \mw \frakg_{-2}& \mw ^{3}\frakg_{-1}& 0 & 0  & 0 
\\ \hline
-4 & \frakg_{-4} & \frakg_{-1} \mw \frakg_{-3} +  \mw ^{2}\frakg_{-2}
&
 \mw ^{2}\frakg_{-1}  \mw  \frakg_{-2} & \mw ^{4}\frakg_{-1} & 0 & 0 
\\ \hline
-5 & 0 & \frakg_{-1} \mw \frakg_{-4} + \frakg_{-2} \mw \frakg_{-3} 
&  \mw ^{2}\frakg_{-1} \mw  \frakg_{-3}
+  \frakg_{-1}  \mw ^{2}\frakg_{-2} &
 \mw ^{3}\frakg_{-1}  \mw   \frakg_{-2} & \mw ^{5}\frakg_{-1} & 0 
\\ \hline
-6 & 0 & \frakg_{-2} \mw \frakg_{-4} +  \mw ^{2}\frakg_{-3} 
& 
(*)
\kmcomment{
\frakg_{-1}  \mw  \frakg_{-2} \mw  \frakg_{-3} 
+  \mw ^{2} \frakg_{-1}  \mw  \frakg_{-4} +  \mw ^{3} \frakg_{-2}
}
& 
(**)
\kmcomment{
 \mw ^{3} \frakg_{-1} \mw  \frakg_{-3} +
 \mw ^{2} \frakg_{-1}  \mw ^{2}\frakg_{-2}
}
&
 \mw ^{4} \frakg_{-1} \mw  \frakg_{-2} 
& 
 \mw ^{6}\frakg_{-1} 
\\ \hline
\end{array}
\]
where \(
(*) = 
\frakg_{-1}  \mw  \frakg_{-2} \mw  \frakg_{-3} 
+  \mw ^{2} \frakg_{-1}  \mw  \frakg_{-4}
+  \mw ^{3} \frakg_{-2}
\)
and \( (**) = 
 \mw ^{3} \frakg_{-1} \mw  \frakg_{-3} +
 \mw ^{2} \frakg_{-1}  \mw ^{2}\frakg_{-2}
\).

The dimensions of those chain spaces are the following: 
\[
\begin{array}{|c||*{7}{c|}}
\hline
w &\dim\CS{w}{1}&\dim\CS{w}{2}&\dim\CS{w}{3}&\dim\CS{w}{4}&\dim\CS{w}{5}
&\dim\CS{w}{6}\\ 
\hline
\hline
-1 & 1 & 0  & 0  & 0  & 0  & 0 
\\ \hline
-2 & 3 & 1  & 0  & 0  & 0  & 0 \\
\hline
-3 & 3 & 3& 1 & 0 & 0  & 0 
\\ \hline
-4 & 1  & 6 & 3 & 1 & 0 & 0 
\\ \hline
-5 & 0 & 10 &  6  & 3 & 1  & 0 
\\ \hline
-6 & 0 & 9 & 11 & 6 & 3 &  1 \\ \hline
\end{array}
\]

\newcommand{\ww}[1]{w_{#1}}

Before discussing individual cases depending on each Lie algebra structure, 
We now study of common behavior of the boundary operator.  
We use impolite notations here like \( \mathbf{1}^{a}\) instead of
\( \mw^{a} \mathbf{1}\) or \( W^{C} =  \ww{1} ^{c_{1}} \ww{2} ^{c_{2}}
\ww{3} ^{c_{3}}  \) instead of \( \mw ^{c_{1}} \ww{1}  \mw ^{c_{2}} \ww{2}
\mw ^{c_{3}} \ww{3} \).  And we abbreviate  
\( Z^{B} = \zb{1}^{b_{1}}  \mw \zb{2}^{j_{2}}  \mw   \zb{3} ^{b_{3}}\) by 
\( Z^{B} = \zb{1}^{b_{1}}  \zb{2}^{j_{2}}   \zb{3} ^{b_{3}}\) where  
\(b_{i}=0,1\). 
Now a typical
generator of the chain space is written as \( \mathbf{1}^{a} \mw Z^{B} \mw
W^{C} \mw V^{\ell} \) or \( \mathbf{1}^{a}  Z^{B}  W^{C}  V^{\ell} \)
where 
\(\ell=0,1\), and whose weight is 
\( a + 2 |B| + 3 |C| + 4 \ell\) and degree is 
\( a +  |B| +  |C| +  \ell\).  
Using properties of \(\pdel\) here,  
\begin{subequations}
\begin{align}
\pdel ( \mathbf{1}^{a} \mw Z^{B} \mw W^{C}) 
&= \pdel( \mathbf{1}^{a}) \mw Z^{B} \mw W^{C}) 
 + (-\mathbf{1})^{a} \pdel ( Z^{B} \mw W^{C}) 
 + \SbtSZ{\mathbf{1}^{a}} {Z^{B} \mw W^{C}} 
 \notag
\\ 
&= 0 + 
  (-\mathbf{1})^{a}( ( \pdel  Z^{B}) \mw W^{C} 
  + \parity{|B|} Z^{B} \mw \pdel  W^{C}  + \SbtSZ{ Z^{B} }{ W^{C}} )
  \notag
  \\&\quad 
 + \SbtSZ{\mathbf{1}^{a}} {Z^{B} \mw W^{C}} 
 \notag
\\ &= (-\mathbf{1})^{a} ( \pdel  Z^{B}) \mw W^{C} 
  + 0 +  0 + \SbtSZ{\mathbf{1}^{a}} {Z^{B} \mw W^{C}} 
  \notag
\\\shortintertext{thus, we have} 
\pdel ( \mathbf{1}^{a} \mw Z^{B} \mw W^{C}) 
&= (-\mathbf{1})^{a} ( \pdel  Z^{B}) \mw W^{C} 
  +  \SbtSZ{\mathbf{1}^{a}} {Z^{B} \mw W^{C}} 
\label{bdary:d3}
\\ 
 \SbtSZ{\mathbf{1}^{a}} {Z^{B} \mw W^{C}} 
&= 
%\textcolor{red}{-} 
 a (-\mathbf{1})^{a-1}  \SbtSZ{\mathbf{1}} {Z^{B} \mw W^{C}} 
\\ 
&=  
%\textcolor{red}{-} 
a (-\mathbf{1})^{a-1} ( 
  \SbtSZ{\mathbf{1}} {Z^{B}} \mw W^{C}  + Z^{B} \mw  
  \SbtSZ{\mathbf{1}} {W^{C}}
  )
\\\shortintertext{and}
\pdel ( \mathbf{1}^{a} \mw Z^{B} \mw W^{C} \mw V^{\ell})&= 
\pdel ( \mathbf{1}^{a} \mw Z^{B} \mw W^{C})\mw V^{\ell} \; .  
\\\shortintertext{And also}
 \SbtSZ{\mathbf{1}} { \zb{i}\mw\zb{j} } &= 
 \SbtSZ{\mathbf{1}} { \zb{i}} \mw\zb{j}  +  \zb{i}\mw
 \SbtSZ{\mathbf{1}}{\zb{j}} \;, \\ 
  \pdel (\zb{i}\mw\zb{j}) & = \SbtD {\zb{i}}{\zb{j}}\;,    
\\
  \SbtSZ{\mathbf{1}} { \zb{1}\mw\zb{2}\mw\zb{3} } &= 
 \SbtSZ{\mathbf{1}} { \zb{1}} \mw\zb{2}\mw\zb{3}  +  
 \zb{1}\mw \SbtSZ{\mathbf{1}}{\zb{2}} \mw\zb{3} + 
 \zb{1}\mw \zb{2}\mw \SbtSZ{\mathbf{1}}{\zb{3}}\;,  
 \\
 \pdel (\zb{1}\mw\zb{2} \mw\zb{3} ) &= \mathop{\mathfrak{S}}_{1,2,3} 
 \SbtD {\zb{1}}{  \zb{2}} \mywedge \zb{3} \;.   
\label{bdary:d3:end}
\end{align}
\end{subequations}
Each subspace  
\( \mw ^{a} \frakg_{-1} \mw^{b}  \frakg_{-2}  \mw ^{c} \frakg_{-3} 
 \mw ^{\ell} \frakg_{-4}  \), 
which is a part of chain
space, and is denoted by  
\( \mathfrak{1}^{a} Z^{b}  W^{c} V ^{\ell} \) in short.   
More precise expression by using basis is given as follows: 
\[ \mathbf{1}^ a Z^{b} W^{c} V^{\ell} = \sum_{|B|=b, |C|=c} 
\lambda_{B,[C,c]}^{a,\ell} 
 \mathbf{1}^ a Z^{B} W^{C} V^{\ell}\] 
where \( B= (b_1,b_2,b_3)\in \{0,1\}^{3}, C= (c_1,c_2,c_3) \in \mN^{3}, 
 |B|=b_1 + b_2 + b_3, 
 |C|=c_1 + c_2 + c_3\), \( \ell=\) 0 or 1,  and \(  
\lambda_{B,[C,c]}^{a,\ell}\) are scalars.  

As described in \eqref{G3:codim:even} or \eqref{G3:codim:odd}, each chain space is given:
%\begin{subequations}
\begin{align*}
\CS{w}{-w-2K} &= 
% \label{G3:codim:even}
  \mathfrak{1} ^{ - w  -3K } Z ^{0} W ^{ K } 
+ \mathfrak{1} ^{ - w -3K - 1 } Z ^{ 2 } W ^{ K  - 1 } 
\notag
%\\&\quad 
+ \mathfrak{1} ^{ - w -3K  } Z ^{ 1 } W ^{ K - 2 } V^{  } 
+ \mathfrak{1} ^{ - w  -3K -1 } Z ^{ 3 } W ^{ K - 3 }V ^{  }
\notag
\\
\CS{w}{-w -2L-1} &= 
% \label{G3:codim:odd}
  \mathfrak{1} ^{-w  -3L -2} Z ^{1} W ^{ L } 
+ \mathfrak{1} ^{-w  -3L - 3 } Z ^{ 3 } W ^{ L  - 1 } 
\notag
%\\&\quad 
+ \mathfrak{1} ^{-w -3L -1 } Z ^{ 0 } W ^{ L - 1 }  V ^{  } 
+ \mathfrak{1} ^{-w  -3L -2 } Z ^{ 2 } W ^{ L - 2 }  V ^{  } 
\; . \notag
\end{align*}
%\end{subequations} 
Thus, using the fact that \( \SbtSZ{\mathbf{1}}{\wb{i}}= \text{u}V\) for
some u,  and so  \( \SbtSZ{\mathbf{1}}{\wb{i}} V = 0\),  the boundary image
is written in general as follows:  
\begin{align}
\shortintertext{When \( -w -m = 2K\), then each element in \(\pdel \CS{w}{m}\) is given by }
X 
=& + \sum\lambda_{[000],[C,K]}^{[-w-3K,0]}( +\SbtSZ{\mathbf{1}^{-w-3K}}{W^{C}})
\label{ukeK:gen} \\&
+\sum\lambda_{[011],[C,K-1]}^{[-w-3K-1,0]}((-\mathbf{1})^{-w-3K-1}\pdel(\zb{2}\zb{3})W^{C}+\SbtSZ{\mathbf{1}^{-w-3K-1}}{\zb{2}\zb{3}W^{C}})
\notag \\&
+\sum\lambda_{[101],[C,K-1]}^{[-w-3K-1,0]}((-\mathbf{1})^{-w-3K-1}\pdel(\zb{1}\zb{3})W^{C}+\SbtSZ{\mathbf{1}^{-w-3K-1}}{\zb{1}\zb{3}W^{C}})
\notag \\&
+\sum\lambda_{[110],[C,K-1]}^{[-w-3K-1,0]}((-\mathbf{1})^{-w-3K-1}\pdel(\zb{1}\zb{2})W^{C}+\SbtSZ{\mathbf{1}^{-w-3K-1}}{\zb{1}\zb{2}W^{C}})
\notag \\&
+\sum\lambda_{[001],[C,K-2]}^{[-w-3K,1]}\SbtSZ{\mathbf{1}^{-w-3K}}{\zb{3}}W^{C}V
\notag \\&
+\sum\lambda_{[010],[C,K-2]}^{[-w-3K,1]}\SbtSZ{\mathbf{1}^{-w-3K}}{\zb{2}}W^{C}V
\notag \\&
+\sum\lambda_{[100],[C,K-2]}^{[-w-3K,1]}\SbtSZ{\mathbf{1}^{-w-3K}}{\zb{1}}W^{C}V
\notag \\&
+\sum\lambda_{[111],[C,K-3]}^{[-w-3K-1,1]}((-\mathbf{1})^{-w-3K-1}
\pdel(\zb{1}\zb{2}\zb{3}) +\SbtSZ{\mathbf{1}^{-w-3K-1}}{\zb{1}\zb{2}
\zb{3}}) W^{C}V
\notag \\\shortintertext{from the original form}
X_{org}&=\sum\lambda_{[000],[C,K]}^{[-w-3K,0]}((-\mathbf{1})^{-w-3K}\pdel(1)W^{C}+\SbtSZ{\mathbf{1}^{-w-3K}}{W^{C}})
\\&
+\sum\lambda_{[011],[C,K-1]}^{[-w-3K-1,0]}((-\mathbf{1})^{-w-3K-1}\pdel(\zb{2}\zb{3})W^{C}+\SbtSZ{\mathbf{1}^{-w-3K-1}}{\zb{2}\zb{3}W^{C}})
\notag \\&
+\sum\lambda_{[101],[C,K-1]}^{[-w-3K-1,0]}((-\mathbf{1})^{-w-3K-1}\pdel(\zb{1}\zb{3})W^{C}+\SbtSZ{\mathbf{1}^{-w-3K-1}}{\zb{1}\zb{3}W^{C}})
\notag \\&
+\sum\lambda_{[110],[C,K-1]}^{[-w-3K-1,0]}((-\mathbf{1})^{-w-3K-1}\pdel(\zb{1}\zb{2})W^{C}+\SbtSZ{\mathbf{1}^{-w-3K-1}}{\zb{1}\zb{2}W^{C}})
\notag \\&
+\sum\lambda_{[001],[C,K-2]}^{[-w-3K,1]}((-\mathbf{1})^{-w-3K}\pdel(\zb{3})W^{C}V+\SbtSZ{\mathbf{1}^{-w-3K}}{\zb{3}W^{C}}V)
\notag \\&
+\sum\lambda_{[010],[C,K-2]}^{[-w-3K,1]}((-\mathbf{1})^{-w-3K}\pdel(\zb{2})W^{C}V+\SbtSZ{\mathbf{1}^{-w-3K}}{\zb{2}W^{C}}V)
\notag \\&
+\sum\lambda_{[100],[C,K-2]}^{[-w-3K,1]}((-\mathbf{1})^{-w-3K}\pdel(\zb{1})W^{C}V+\SbtSZ{\mathbf{1}^{-w-3K}}{\zb{1}W^{C}}V)
\notag \\&
+\sum\lambda_{[111],[C,K-3]}^{[-w-3K-1,1]}((-\mathbf{1})^{-w-3K-1}\pdel(\zb{1}\zb{2}\zb{3})W^{C}V+\SbtSZ{\mathbf{1}^{-w-3K-1}}{\zb{1}\zb{2}\zb{3}W^{C}}V)
\notag
\end{align}
\begin{align}
\shortintertext{Similarly, when \( -w -m = 2L+1\), then each element in \(\pdel \CS{w}{m}\) is given by }
%ukeAllLg 
Y= &+ \sum\lambda_{[001],[C,L]}^{[-w-3L-2,0]} \SbtSZ{\mathbf{1}^{-w-3L-2}}{\zb{3}W^{C}}
\label{ukeL:gen}
\\&
+\sum\lambda_{[010],[C,L]}^{[-w-3L-2,0]} \SbtSZ{\mathbf{1}^{-w-3L-2}}{\zb{2}W^{C}}
\notag \\&
+\sum\lambda_{[100],[C,L]}^{[-w-3L-2,0]} \SbtSZ{\mathbf{1}^{-w-3L-2}}{\zb{1}W^{C}}
\notag \\&
+\sum\lambda_{[111],[C,L-1]}^{[-w-3L-3,0]}(
(-\mathbf{1})^{-w-3L-3}\pdel(\zb{1}\zb{2}\zb{3})W^{C}
+\SbtSZ{\mathbf{1}^{-w-3L-3}}{\zb{1}\zb{2}\zb{3}W^{C}})
\notag \\&
\kmcomment{
+\sum\lambda_{[000],[C,L-1]}^{[-w-3L-1,1]}(\SbtSZ{\mathbf{1}^{-w-3L-1}}{W^{C}}V)
\notag \\&
}
+\sum\lambda_{[011],[C,L-2]}^{[-w-3L-2,1]}((-\mathbf{1})^{-w-3L-2}
\pdel(\zb{2}\zb{3}) 
+\SbtSZ{\mathbf{1}^{-w-3L-2}}{\zb{2}\zb{3}}) W^{C}V
\notag \\&
+\sum\lambda_{[101],[C,L-2]}^{[-w-3L-2,1]}((-\mathbf{1})^{-w-3L-2}
\pdel(\zb{1}\zb{3})
+\SbtSZ{\mathbf{1}^{-w-3L-2}}{\zb{1}\zb{3}}) W^{C}V
\notag \\&
+\sum\lambda_{[110],[C,L-2]}^{[-w-3L-2,1]}((-\mathbf{1})^{-w-3L-2}
\pdel(\zb{1}\zb{2}) 
+\SbtSZ{\mathbf{1}^{-w-3L-2}}{\zb{1}\zb{2}} )W^{C}V
\notag
\end{align}

\begin{comment}
\begin{align*}
\shortintertext{Similarly, when \( -w -m = 2L+1\), then each element in \(\pdel \CS{w}{m}\) is given by }
ukeAllLg &=\sum\lambda_{[001],[C,L]}^{[-w-3L-2,0]}((-\mathbf{1})^{-w-3L-2}\pdel(\zb{3})W^{C}+\SbtSZ{\mathbf{1}^{-w-3L-2}}{\zb{3}W^{C}})
\\&
+\sum\lambda_{[010],[C,L]}^{[-w-3L-2,0]}((-\mathbf{1})^{-w-3L-2}\pdel(\zb{2})W^{C}+\SbtSZ{\mathbf{1}^{-w-3L-2}}{\zb{2}W^{C}})
\\&
+\sum\lambda_{[100],[C,L]}^{[-w-3L-2,0]}((-\mathbf{1})^{-w-3L-2}\pdel(\zb{1})W^{C}+\SbtSZ{\mathbf{1}^{-w-3L-2}}{\zb{1}W^{C}})
\\&
+\sum\lambda_{[111],[C,L-1]}^{[-w-3L-3,0]}((-\mathbf{1})^{-w-3L-3}\pdel(\zb{1}\zb{2}\zb{3})W^{C}+\SbtSZ{\mathbf{1}^{-w-3L-3}}{\zb{1}\zb{2}\zb{3}W^{C}})
\\&
+\sum\lambda_{[000],[C,L-1]}^{[-w-3L-1,1]}((-\mathbf{1})^{-w-3L-1}\pdel(1)W^{C}V+\SbtSZ{\mathbf{1}^{-w-3L-1}}{W^{C}}V)
\\&
+\sum\lambda_{[011],[C,L-2]}^{[-w-3L-2,1]}((-\mathbf{1})^{-w-3L-2}\pdel(\zb{2}\zb{3})W^{C}V+\SbtSZ{\mathbf{1}^{-w-3L-2}}{\zb{2}\zb{3}W^{C}}V)
\\&
+\sum\lambda_{[101],[C,L-2]}^{[-w-3L-2,1]}((-\mathbf{1})^{-w-3L-2}\pdel(\zb{1}\zb{3})W^{C}V+\SbtSZ{\mathbf{1}^{-w-3L-2}}{\zb{1}\zb{3}W^{C}}V)
\\&
+\sum\lambda_{[110],[C,L-2]}^{[-w-3L-2,1]}((-\mathbf{1})^{-w-3L-2}\pdel(\zb{1}\zb{2})W^{C}V+\SbtSZ{\mathbf{1}^{-w-3L-2}}{\zb{1}\zb{2}W^{C}}V);
\end{align*}
\end{comment}

Since both \(\pdel\) and \( \SbtSZ{\cdot}{\cdot}\)  depend on the Lie
algebra structure, our discussion needs to be developed individually.  
\kmcomment{
We use the notation \(\zb{i}\) instead of \(\sigma^{i}\) in
\eqref{ext:deriv:d} and denote the 2-form \( \zb{i} \wedge \zb{i+1}\) by \(
\wb{i+2}\), where indices are reduced by modulo 3.
Put \( \zb{1} \wedge \zb{2} \wedge \zb{3}\) by \(V\). 
}

\subsubsection{Case of SO(3):} 

\begin{wrapfigure}[4]{r}[5mm]{0.45\textwidth}
\vspace{-1.0 \baselineskip}
%\begin{center}
\( \renewcommand{\arraystretch}{0.8}
\begin{array}{c|c|*{3}{c}|*{3}{c}|c}
\hline
& \mathbf{1} &
\zb{1} &
\zb{2} &
\zb{3} &
\wb{1} &
\wb{2} &
\wb{3} &
V
\\\hline
\mathbf{1} & 0 & -2 \wb{1} & -2\wb{2} & -2\wb{3} & 0 & 0 & 0 &  \\
\cline{1-8}
\zb{1} & 2 \wb{1} & 0 & 0 & 0 & \\
\zb{2} & 2 \wb{2} & 0 & 0 & 0 & \\
\zb{3} & 2 \wb{3} & 0 & 0 & 0 &  \\
\cline{1-5}
\wb{1} & 0&  \\
\wb{2} & 0 & \\
\wb{3} & 0 &  \\
\cline{1-2}
V       
\end{array}
\)
\end{wrapfigure}
Now we have \[ d \zb{i} =
- 2 \wb{i} \quad i=1,2,3\] and the brackets of superalgebra are given as
follows:
%\begin{center}

\(
\begin{array}{*{5}{l}}
\SbtD{1}{1} = 0 & \SbtD{1}{\zb{j}} = - 2 \ww{j} & \SbtD{1}{\ww{j}} = 0 & 
\SbtD{1}{V} = 0 \\ 
\SbtD{\zb{i}}{1} = 2\ww{i} & \SbtD{\zb{i}}{\zb{j}} = 0 & \SbtD{\zb{i}}{\ww{j}} = 0 & 
\SbtD{\zb{i}}{V} = 0 \\ 
\SbtD{\ww{i}}{1} = 0 & \SbtD{\ww{i}}{\zb{j}} = 0 & \SbtD{\ww{i}}{\ww{j}} = 0 & 
\SbtD{\ww{i}}{V} = 0 \\ 
\SbtD{V}{1} = 0 & \SbtD{V}{\zb{j}} = 0 & \SbtD{V}{\ww{j}} = 0 & 
\SbtD{V}{V} = 0 \\ 
\end{array}
\)
%\end{center}

\( V\) and \(\ww{i}\) are central elements of this algebra.  
\eqref{bdary:d3}\(\sim\)
\eqref{bdary:d3:end}
imply 
\begin{align*}
\pdel(   \mathbf{1} ^{a}  Z^{B}   W^{C} V ^{\ell} )
&= 
 a ( - \mathbf{1} )^{a-1} \SbtSZ { \mathbf{1} }{ Z^{B}}  W^{C}  V  ^{\ell}  
\\
\SbtSZ{ \mathbf{1}} { \zb{j} } &=  -2 \ww{j}
\\
\SbtSZ{ \mathbf{1}} { \zb{j}  \mw  \zb{k}  } &= 
 ( \zb{j}  \mw   ( -2 \ww{k} )
+  \zb{k}  \mw  ( -2 \ww{j} )) 
\\
\SbtSZ{ \mathbf{1}} {
\zb{1}  \mw   \zb{2}  \mw   \zb{3} }
&= 
 \SbtD{ \mathbf{1} }{\zb{1} } \mw   \zb{2}  \mw   \zb{3} +
\zb{1}  \mw  \SbtD{\mathbf{1}}{ \zb{2}}  \mw   \zb{3} +
\zb{1}  \mw   \zb{2}  \mw  \SbtD{\mathbf{1}}{ \zb{3} }
\\&= 
 \zb{1}  \mw   \zb{2}  \mw  ( -2 \ww{3} )
- \zb{1}  \mw   \zb{3}  \mw  ( -2 \ww{2} )  
+ \zb{2}  \mw   \zb{3}  \mw  ( -2 \ww{1} )
\end{align*}
%Discussion of rank here: 
We remember the decomposition of chain spaces by \eqref{G3:codim:even} and
\eqref{G3:codim:odd}, and knowing the properties of the boundary operator
which is defined by the Lie algebra structure, we first have direct sum
decomposition of subspaces as below: 
%%%%%%%%%%%%%%%%%%%%%%%%%%%%%%%%%%%%%%%%%%%%%%%%%%%%%
\begin{align*}
\shortintertext{When \( m= -w -2K\), a general element in 
\(\pdel \CS{w}{m}\) is written as}  
X =&
+2 \sum\lambda_{[011],[C,K-1]}^{[-w-3K-1,0]}(-w-3K-1)(-II)^{-w-3K-2}(\zb{3}\wb{2}-\zb{2}\wb{3})W^{C}
\\&
+2\sum\lambda_{[101],[C,K-1]}^{[-w-3K-1,0]}(-w-3K-1)(-II)^{-w-3K-2}(\zb{3}\wb{1}-\zb{1}\wb{3})W^{C}
\\&
+2\sum\lambda_{[110],[C,K-1]}^{[-w-3K-1,0]}(-w-3K-1)(-II)^{-w-3K-2}(\zb{2}\wb{1}-\zb{1}\wb{2})W^{C}
\\&
-2\sum\lambda_{[001],[C,K-2]}^{[-w-3K,1]}(-w-3K)(-II)^{-w-3K-1}\wb{3}W^{C}V
\\&
-2\sum\lambda_{[010],[C,K-2]}^{[-w-3K,1]}(-w-3K)(-II)^{-w-3K-1}\wb{2}W^{C}V
\\&
-2\sum\lambda_{[100],[C,K-2]}^{[-w-3K,1]}(-w-3K)(-II)^{-w-3K-1}\wb{1}W^{C}V
\\&
-2\sum\lambda_{[111],[C,K-3]}^{[-w-3K-1,1]}(-w-3K-1)(-II)^{-w-3K-2}(\zb{2}\zb{3}\wb{1}-\zb{1}\zb{3}\wb{2}+\zb{1}\zb{2}\wb{3})W^{C}V
\\
\shortintertext{And when \( m= -w -2L-1\), a general element in 
\(\pdel \CS{w}{m}\) is written as}  
Y =& -2 \sum\lambda_{[001],[C,L]}^{[-w-3L-2,0]}(-w-3L-2)(-II)^{-w-3L-3}\wb{3} W^{C}
\\&
- 2 \sum\lambda_{[010],[C,L]}^{[-w-3L-2,0]}(-w-3L-2)(-II)^{-w-3L-3}\wb{2}W^{C}
\\&
-2\sum\lambda_{[100],[C,L]}^{[-w-3L-2,0]}(-w-3L-2)(-II)^{-w-3L-3} \wb{1} W^{C}
\\&
-2\sum\lambda_{[111],[C,L-1]}^{[-w-3L-3,0]}(-w-3L-3)(-II)^{-w-3L-4}(\zb{2}\zb{3} \wb{1} -\zb{1}\zb{3} \wb{2} +\zb{1}\zb{2}\wb{3})W^{C}
\\&
-2 \sum\lambda_{[011],[C,L-2]}^{[-w-3L-2,1]}(-w-3L-2)(-II)^{-w-3L-3}
(-\zb{3}\wb{2}+\zb{2} \wb{3} )W^{C}V
\\&
+ 2 \sum\lambda_{[101],[C,L-2]}^{[-w-3L-2,1]}(-w-3L-2)(-II)^{-w-3L-3}
(\zb{3}\wb{1}-\zb{1}\wb{3})W^{C}V
\\&
+2 \sum\lambda_{[110],[C,L-2]}^{[-w-3L-2,1]}(-w-3L-2)(-II)^{-w-3L-3}
( \zb{2}\wb{1}- \zb{1}\wb{2})W^{C}V
\end{align*}
%
%%%%%%%%%%%%%%%%%%%%%%%%%%%%%%%%%%%%%%%%%%%%%%%%%%%%%
\begin{align*}
\shortintertext{Thus, when \( m= -w -2K\)}
\dim \pdel \CS{w}{m} &= 
3 \tbinom{-w-3K-2}{0} \tbinom{2+K-1}{K-1}
+ 3 \tbinom{-w-3K-1}{0}  \tbinom{2+K-2}{K-2} 
+ \tbinom{-w-3K-2}{0}  \tbinom{2+K-3}{K-3} \;, 
\\
\shortintertext{and when \( m= -w -2L-1\)}
\dim \pdel \CS{w}{m} &= 
3 \tbinom{-w-3L-3}{0} \tbinom{2+L}{L}
+ \tbinom{-w-3L-4}{0}  \tbinom{2+L-1}{L-1} 
+ 3 \tbinom{-w-3L-3}{0}  \tbinom{2+L-2}{L-2} \;. 
\end{align*}

%%%%%%%%%%%%%%%%%%%%%%%%%%%%%%%%%%%%%%%%%%%%%%%%%%%%%
%%%%%%%%%%%%%%%%%%%%%%%%%%%%%%%%%%%%%%%%%%%%%%%%%%%%%
\subsubsection{\(SL(2,\mR)\):} 
In the case of \(SL(2,\mR)\), 
we have \( d \zb{1} =  2 \ww{1}\;,\;  
d \zb{i} = - 2 \ww{i} \quad i=1,2,3\) and 
the brackets of superalgebra are given as follows:
\[
\begin{array}{*{5}{l}}
\SbtD{1}{1} = 0 & \SbtD{1}{\zb{j}} = -2(1 - 2\delta_{1j})\ww{j} & \SbtD{1}{\ww{j}} = 0 & 
\SbtD{1}{V} = 0 \\ 
\SbtD{\zb{i}}{1} = 2(1-2\delta_{1i} )\ww{i} & \SbtD{\zb{i}}{\zb{j}} = 0 & \SbtD{\zb{i}}{\ww{j}} = 0 & 
\SbtD{\zb{i}}{V} = 0 \\ 
\SbtD{\ww{i}}{1} = 0 & \SbtD{\ww{i}}{\zb{j}} = 0 & \SbtD{\ww{i}}{\ww{j}} = 0 & 
\SbtD{\ww{i}}{V} = 0 \\ 
\SbtD{V}{1} = 0 & \SbtD{V}{\zb{j}} = 0 & \SbtD{V}{\ww{j}} = 0 & 
\SbtD{V}{V} = 0 \\ 
\end{array}
\]
\( V\) and \(\ww{i}\) are central elements of this algebra. 
We may say this algebra is isomorphic to 
the algebra associated with \(SO(3,\mR)\).

\medskip

\subsubsection{\(\dim \Sbt{\frakg}{\frakg} = 2\):}

\begin{wrapfigure}[7]{r}[5mm]{0.65 \textwidth}
\vspace{-2.8 \baselineskip}
%\begin{center}

%\begin{center}
\(\renewcommand{\arraystretch}{0.8}
\begin{array}{c|c|*{3}{c}|*{3}{c}|c}
\hline
& \mathbf{1} &
\zb{1} &
\zb{2} &
\zb{3} &
\wb{1} &
\wb{2} &
\wb{3} &
V
\\\hline
\mathbf{1} & 0 & 2 \wb{2} & -2\kappa\wb{1} & 0 & 0 & 0 & 2(1+\kappa)V &  \\
\cline{1-8}
\zb{1} & -2 \wb{2} & 0 & -2(1+\kappa)V & 0 & \\
\zb{2} & 2\kappa \wb{1} & 2(1+\kappa)V & 0 & 0 & \\
\zb{3} & 0 & 0 & 0 & 0 &  \\
\cline{1-5}
\wb{1} & 0&  \\
\wb{2} & 0 & \\
\wb{3} & 2(1+\kappa)V &  \\
\cline{1-2}
V       
\end{array}
\)
%\end{center}
\end{wrapfigure}
In the case of \(\dim \Sbt{\frakg}{\frakg} = 2\), we may take a basis 
\( \Sbt{\xi_{1}}{\xi_{2}} = 0\;,\; 
\Sbt{\xi_{2}}{\xi_{3}} = \kappa \xi_{2}\;,\; 
\Sbt{\xi_{3}}{\xi_{1}} = - \xi_{1}\) where \(\kappa \ne 0\), 
and we have \( d \zb{1} =  2 \ww{2}\;,\;  d \zb{2} =  -2 \kappa \ww{1}\;,\;
d \zb{3} = 0 \).  The brackets of superalgebra are given as follows:

\eqref{bdary:d3} implies 
\begin{align*}
\pdel ( \mathbf{1}^{a}  Z^{B} W^{C} V^{\ell} )
&= 
 (-\mathbf{1})^{a} ( \pdel Z^{B}) W^{C} V^{\ell}
 + a (-\mathbf{1})^{a-1}(
 \SbtSZ{ \mathbf{1} }{ Z^{B} } W^{C} + Z^{B} 
 \SbtSZ{ \mathbf{1} }{ W^{C} } ) V^{\ell}
\end{align*}
because \( \wb{1},\wb{2}\) are central.  
From now on, we abbreviate \( \zb{i} \mw \zb{j}\) as \( \zb{i,j}\)
and \( \zb{1} \mw \zb{2} \mw \zb{3}\) as \( \zb{1,2,3}\). 
\(\pdel Z^{B} = 0\) if \(|B| \leqq 1\), 
\( \pdel ( \zb{i,j}) =
\SbtD{ \zb{i} }{\zb{j}} = \begin{cases} -2 (1+\kappa) V & (i,j) = (1,2) \\ 0
& \text{otherwise} \end{cases}\), \( \pdel ( \zb{1,2,3} )
= 2 (1+\kappa) \zb{3} \mw  V\).  \( \SbtSZ{ \mathbf{1}^a }{Z^B} = a ( -
\mathbf{1} )^{a-1} \SbtSZ{ \mathbf{1} }{Z^B} \).   
\kmcomment{
\( \pdel ( \zb{i}\mw \zb{j}) =
\SbtD{ \zb{i} }{\zb{j}} = \begin{cases} -2 (1+\kappa) V & (i,j) = (1,2) \\ 0
& \text{otherwise} \end{cases}\), \( \pdel ( \zb{1} \mw  \zb{2}  \mw  \zb{3})
= 2 (1+\kappa) \zb{3} \mw  V\).  \( \SbtSZ{ \mathbf{1}^a }{Z^B} = a ( -
\mathbf{1} )^{a-1} \SbtSZ{ \mathbf{1} }{Z^B} \).   
}
\(\SbtSZ{\mathbf{1}}{ \zb{i,j} } = \SbtSZ{\mathbf{1}}{ \zb{i}}
\mw  \zb{j} + \zb{i}  \mw  \SbtSZ{\mathbf{1}}{ \zb{j}  } \) and
\(\SbtSZ{\mathbf{1}}{ \zb{1,2,3}   } = 
\kmcomment{
\SbtSZ{\mathbf{1}}{ \zb{1}} \mw \zb{2} \mw \zb{3} + \zb{1} \mw
\SbtSZ{\mathbf{1}}{ \zb{2}} \mw \zb{3} + \zb{1} \mw \zb{2} \mw
\SbtSZ{\mathbf{1}}{ \zb{3}} =}
2  \zb{2} \mw \zb{3} \mw \wb{2}
+ 2\kappa \zb{1} \mw \zb{3} \mw \wb{1}
\), and  
\kmcomment{
Since \( \wb{1},\wb{2}\) are central, 
\begin{align*}
\SbtSZ { \mathbf{1}^{a}  Z^{B}}{W^{C}} 
& = c_{3} \SbtSZ { \mathbf{1}^{a}  Z^{B}}{\wb{3}}\mw  W^{c_1,c_2,c_3 -1} 
 = c_{3} a (-\mathbf{1})^{a-1} Z^{B} 
 \SbtSZ { \mathbf{1}}{\wb{3}}\mw  W^{c_1,c_2,c_3 -1} 
 \\& 
 = 2  c_{3} (1+\kappa) a (-\mathbf{1})^{a-1} Z^{B} \mw
  V \mw  W^{c_1,c_2,c_3 -1}  \; . 
\end{align*}
}
\[ \SbtSZ { \mathbf{1}}{W^{C}}
 = c_{3} 
 \SbtSZ { \mathbf{1}}{\wb{3}}  W^{c_1,c_2,c_3 -1} 
 = 2  c_{3} (1+\kappa)  V \mw  W^{c_1,c_2,c_3 -1}  
 = 2  (1+\kappa)  V \mw \frac{c_{3}}{\wb{3}} W^{C}  \;. \] 
Using
\eqref{ukeK:gen}, we have an expression of  
\(\pdel\)-image  for each generator in \eqref{G3:codim:even} as follows:  
\begin{subequations}
\begin{align}
\shortintertext{When \( -w -m = 2K\), }
Kd=&
(-w-3K)
(-\mathbf{1})^{-w-3K-1}\sum\lambda_{[000],[C,K]}^{[-w-3K,0]}
(\SbtSZ{\mathbf{1}}{W^{C}})
\label{kk:gen:a}
\\& +
(-w-3K-1)
(-\mathbf{1})^{-w-3K-2}\sum\lambda_{[011],[C,K-1]}^{[-w-3K-1,0]}
(2\kappa{}\zb{3}\wb{1}W^{C}+\zb{2}\zb{3}\SbtSZ{\mathbf{1}}{W^{C}})
\label{kk:gen:b}
\\& +
(-w-3K-1)
(-\mathbf{1})^{-w-3K-2}\sum\lambda_{[101],[C,K-1]}^{[-w-3K-1,0]}
(-2\zb{3}\wb{2}W^{C}+\zb{1}\zb{3}\SbtSZ{\mathbf{1}}{W^{C}})
\label{kk:gen:c}
\\& +(-\mathbf{1})^{-w-3K-2}\sum\lambda_{[110],[C,K-1]}^{[-w-3K-1,0]}
(\mathbf{1}(1+\kappa{}  ) 2 V )W^{C}
\label{kk:gen:d}
\\& \qquad \qquad \qquad 
+(-w-3K-1)((-2\kappa{}\wb{1}\zb{1}-2\wb{2}\zb{2})W^{C}+\zb{1}\zb{2}\SbtSZ{\mathbf{1}}{W^{C}}))
\notag 
%\\& +(-\mathbf{1})^{-w-3K-1}\sum\lambda_{[001],[C,K-2]}^{[-w-3K,1]} (0 ) 
%\notag
\\& +
(-w-3K)
(-\mathbf{1})^{-w-3K-1}\sum\lambda_{[010],[C,K-2]}^{[-w-3K,1]}
(-2\kappa{}\wb{1}W^{C}  V)
\label{kk:gen:e}
\\& +
(-w-3K)
(-\mathbf{1})^{-w-3K-1}\sum\lambda_{[100],[C,K-2]}^{[-w-3K,1]}
( 2\wb{2}W^{C} V)
\label{kk:gen:f}
\\& +
(-w-3K-1)
(-\mathbf{1})^{-w-3K-2}\sum\lambda_{[111],[C,K-3]}^{[-w-3K-1,1]}
( (2\kappa{}\wb{1}\zb{1}\zb{3}+2\wb{2}\zb{2}\zb{3})W^{C} V)
\label{kk:gen:g}
\end{align}
\end{subequations}

Assume \( -w -m = 2K\) and \(\kappa+1=0\). Then \(\SbtSZ{1}{W^C}=0\) and the
expression above becomes simpler like 
\begin{align*}
Kdz = &
(-w-3K-1)
(-\mathbf{1})^{-w-3K-2}\sum\lambda_{[011],[C,K-1]}^{[-w-3K-1,0]}
(2\kappa{}\zb{3}\wb{1}W^{C} 
)
\\& +
(-w-3K-1)
(-\mathbf{1})^{-w-3K-2}\sum\lambda_{[101],[C,K-1]}^{[-w-3K-1,0]}
(-2\zb{3}\wb{2}W^{C}
)
\\& +(-\mathbf{1})^{-w-3K-2}\sum\lambda_{[110],[C,K-1]}^{[-w-3K-1,0]}
(-w-3K-1)(-2\kappa{}\wb{1}\zb{1}-2\wb{2}\zb{2})W^{C}
%\\& +(-\mathbf{1})^{-w-3K-1}\sum\lambda_{[001],[C,K-2]}^{[-w-3K,1]} (0 ) 
%\notag
\\& +
(-w-3K)
(-\mathbf{1})^{-w-3K-1}\sum\lambda_{[010],[C,K-2]}^{[-w-3K,1]}
(-2\kappa{}\wb{1}W^{C}  V)
\\& +
(-w-3K)
(-\mathbf{1})^{-w-3K-1}\sum\lambda_{[100],[C,K-2]}^{[-w-3K,1]}
( 2\wb{2}W^{C} V)
\\& +
(-w-3K-1)
(-\mathbf{1})^{-w-3K-2}\sum\lambda_{[111],[C,K-3]}^{[-w-3K-1,1]}
 (2\kappa{}\wb{1}\zb{1}\zb{3}+2\wb{2}\zb{2}\zb{3})W^{C} V
\end{align*}
If \( -w -3K -1 <0\), then \(Kdz =0\), i.e., \( \dim \CS{w}{m} = 0\). 
If \( -w -3K -1 =0\), then \[Kdz =
(-w-3K) (-\mathbf{1})^{-w-3K-1} ( 
\sum\lambda_{[010],[C,K-2]}^{[-w-3K,1]}
(-2\kappa{}\wb{1}W^{C}  V)
 + \sum\lambda_{[100],[C,K-2]}^{[-w-3K,1]}
( 2\wb{2}W^{C} V)
\] i.e., the space is a union of two subspaces of dimension  
\(  \tbinom{2+K-2}{K-2}\) and those intersection is \( \tbinom{2+K-3}{K-3} \)
-dimensional. Thus,  
\( \dim \CS{w}{m} = 2 \tbinom{2+K-2}{K-2}- \tbinom{2+K-3}{K-3} \). 

If \( -w -3K -1 >0\), then
we apply the same discussion for the first two subspaces, and we get 
\begin{align*}
 \dim \CS{w}{m} = &  
2 \tbinom{2+K-1}{K-1} - \tbinom{2+K-2}{K-2} 
 + \tbinom{2+K-1}{K-1} 
+ 2 \tbinom{2+K-2}{K-2} - \tbinom{2+K-3}{K-3} 
 + \tbinom{2+K-3}{K-3} 
 \\ = & 
3 \tbinom{2+K-1}{K-1} + \tbinom{2+K-2}{K-2} 
\end{align*}
Assume \( -w -m = 2K\) and \(\kappa+1 \ne 0\). We use the formula 
\eqref{kk:gen:a} \(\sim\) \eqref{kk:gen:g}.  

If \( -w -3K-1 < 0\) then the source element is zero and Kd = 0. 
If \( -w -3K-1 = 0\) then we have 
\begin{align*}
X =&
\sum\lambda_{[000],[C,K]}^{[-w-3K,0]}\SbtSZ{\mathbf{1}}{W^{C}}
 - \sum\lambda_{[110],[C,K-1]}^{[-w-3K-1,0]} 2(\kappa{}  + 1 )  V W^{C}
\\& 
+2 ( -  \sum\lambda_{[010],[C,K-2]}^{[-w-3K,1]}\kappa{}\wb{1}W^{C}
 + \sum\lambda_{[100],[C,K-2]}^{[-w-3K,1]}\wb{2}W^{C})V
 \\=& 2
\sum ( \lambda_{[000],[[c_1,c_2,1+c_3],K]}^{[-w-3K,0]} (1+\kappa) (1+c_3)
 + \lambda_{[110],[[c_1,c_2,c_3],K-1]}^{[-w-3K-1,0]} (1 + \kappa{})
\\& 
-  \lambda_{[010],[[-1+c_1, c_2, c_3],K-2]}^{[-w-3K,1]}\kappa{}
 + \sum\lambda_{[100],[[c_1,-1+c_2,c_3],K-2]}^{[-w-3K,1]})
 W^{c_1,c_2,c_3} V
 \\\shortintertext{
Thus, the maximal linearly independent generators are } 
& \lambda_{[000],[[c_1,c_2,1+c_3],K]}^{[-w-3K,0]} (1+\kappa) (1+c_3)
 + \lambda_{[110],[[c_1,c_2,c_3],K-1]}^{[-w-3K-1,0]} (1 + \kappa{})
 \\& 
-  \lambda_{[010],[[-1+c_1, c_2, c_3],K-2]}^{[-w-3K,1]}\kappa{}
 + \sum\lambda_{[100],[[c_1,-1+c_2,c_3],K-2]}^{[-w-3K,1]}
\end{align*}
with \(c_1 + c_2 +c_3 = K-1\), and so the rank is \( \tbinom{2+K-1}{K-1}
\).   

If \( -w -3K-1 > 0\) then the above shows they are divided into 4 linearly
independent parts: The first consists of \eqref{kk:gen:a}, \eqref{kk:gen:e}
and \eqref{kk:gen:f}, which
have the common factor \(V\), the second is  \eqref{kk:gen:b} and   
\eqref{kk:gen:c}  which have the
common factor \(\zb{3} \), the third part is \eqref{kk:gen:d} whose rank is
\(\tbinom{2+K-1}{K-1}\), 
the fourth part is \eqref{kk:gen:g} whose rank is \(\tbinom{2+K-3}{K-3}\). 
Since 
\begin{align*} 
\eqref{kk:gen:a} + \eqref{kk:gen:e} + \eqref{kk:gen:f} = & 
- (-w-3K) (-\mathbf{1})^{-w-3K-1}2 (  \parity{K}\sum\lambda_{[000],[C,K]}^{[-w-3K,0]}
 ( 1+\kappa ) \frac{c_3}{\wb{3}} W^{C}
\\& 
 +  \sum\lambda_{[010],[C,K-2]}^{[-w-3K,1]}\kappa{}\wb{1}W^{C}
 - \sum\lambda_{[100],[C,K-2]}^{[-w-3K,1]}\wb{2}W^{C}) ) V
 \\=& 
- (-w-3K) (-\mathbf{1})^{-w-3K-1}2 \sum_{|C|=K-1} (  \parity{K}\lambda_{[000],[[c_1,
c_2,1+c_3],K]}^{[-w-3K,0]}
 ( 1+\kappa )(1+c_3)
\\& 
 + \lambda_{[010],[[-1+c_1,c_2,c_3],K-2]}^{[-w-3K,1]}\kappa{}
 - \lambda_{[100],[[c_1, -1+c_2, c_3],K-2]}^{[-w-3K,1]}) W^{c_1,c_2,c_3}  V
 \; , 
\end{align*}
the rank of the first part
\( \eqref{kk:gen:a} + \eqref{kk:gen:e} + \eqref{kk:gen:f} \) 
is \(\tbinom{2+K-1} {K-1}\). 
We get 
\begin{align*} 
\eqref{kk:gen:b} + \eqref{kk:gen:c}  = & 
 +
(-w-3K-1)
(-\mathbf{1})^{-w-3K-2}\zb{3} \sum\lambda_{[011],[C,K-1]}^{[-w-3K-1,0]}
(2\kappa{}\wb{1}W^{C}-\zb{2}\SbtSZ{\mathbf{1}}{W^{C}})
%\label{kk:gen:b}
\\& +
(-w-3K-1)
(-\mathbf{1})^{-w-3K-2} \zb{3} \sum\lambda_{[101],[C,K-1]}^{[-w-3K-1,0]}
(-2\wb{2}W^{C}-\zb{1}\SbtSZ{\mathbf{1}}{W^{C}})
%\label{kk:gen:c}
\\= &
 (-w-3K-1) (-\mathbf{1})^{-w-3K-2} 2 \zb{3} \mw 
\\&\qquad 
(  
 \sum\lambda_{[011],[C,K-1]}^{[-w-3K-1,0]}
( \kappa{}\wb{1}W^{C}-\zb{2}   
       (1+\kappa)  V  \frac{c_{3}}{\wb{3}} W^{C}  
)
\\& \qquad\qquad  - \sum\lambda_{[101],[C,K-1]}^{[-w-3K-1,0]}
( \wb{2}W^{C}+\zb{1}     
     (1+\kappa)  V  \frac{c_{3}}{\wb{3}} W^{C}  
)
 \; , 
\end{align*}
We claim that the rank of 
\eqref{kk:gen:b} + \eqref{kk:gen:c}  is \( 2 \tbinom{2+K-1}{K-1} -
\tbinom{1+K-2}{K-2}\).  Reason is as follows. We have two kinds of
generators:  
\( 
( \kappa{}\wb{1}W^{C}-\zb{2} (1+\kappa)  V  \frac{c_{3}}{\wb{3}} W^{C}  )\) 
and \( ( \wb{2}W^{C'}+\zb{1} (1+\kappa)  V  \frac{c'_{3}}{\wb{3}} W^{C'}  )
\) where \( |C| = |C'| = K-1\). Since \( 1+ \kappa \ne 0\), if \( c_3 \ne 0\)
or \( c'_3 \ne 0\) then they are linearly independent, and there  
the rank is \(
2 ( \tbinom{2+K-1}{K-1} - \tbinom{1+K-1}{K-1} \). 
Otherwise, we can
recall the same old discussion and the rank is  \( 
2 \tbinom{1+K-1}{K-1} - \tbinom{1+K-2}{K-2} \). 
Finally, \( \dim \pdel \CS{w}{m}= 
4 \tbinom{2+K-1}{K-1} - \tbinom{1+K-2}{K-2} 
+ \tbinom{2+K-3}{K-3}  \) when \( -w -m = 2K\) and \( 1+\kappa \ne 0\).

We study when \( -w -m = 2L + 1\). Since \( \SbtSZ{\mathbf{1}}{\zb{i}} V = 0
\), we have  
\begin{align*}
Ld &=
(-w-3L-2) (-\mathbf{1})^{-w-3L-3}\sum\lambda_{[001],[C,L]}^{[-w-3L-2,0]}
(\zb{3}\SbtSZ{\mathbf{1}}{W^{C}})
\\& +(-w-3L-2)(-\mathbf{1})^{-w-3L-3}\sum\lambda_{[010],[C,L]}^{[-w-3L-2,0]}
((-2\kappa{}\wb{1}W^{C}+\zb{2}\SbtSZ{\mathbf{1}}{W^{C}}))
\\& +(-w-3L-2)(-\mathbf{1})^{-w-3L-3}\sum\lambda_{[100],[C,L]}^{[-w-3L-2,0]}
((2\wb{2}W^{C}+\zb{1}\SbtSZ{\mathbf{1}}{W^{C}}))
\\& +
\sum\lambda_{[111],[C,L-1]}^{[-w-3L-3,0]}
(
(-\mathbf{1})^{-w-3L-3}
(2(1+\kappa)\zb{3} V )W^{C}
\\&\qquad 
+(-w-3L-3)(-\mathbf{1})^{-w-3L-4}
((2\kappa{}\wb{1}\zb{1}\zb{3}+2\wb{2}\zb{2}\zb{3})W^{C}+\zb{1}\zb{2}\zb{3}\SbtSZ{\mathbf{1}}{W^{C}}))
%% \\& +(-\mathbf{1})^{-w-3L-2}\sum\lambda_{[000],[C,L-1]}^{[-w-3L-1,1]}
%% ((-w-3L-1)\SbtSZ{\mathbf{1}}{W^{C}}V)
\\& +(-w-3L-2)(-\mathbf{1})^{-w-3L-3}\sum\lambda_{[011],[C,L-2]}^{[-w-3L-2,1]}
((2\zb{3}\kappa{}\wb{1}W^{C})V)
\\& +(-w-3L-2)(-\mathbf{1})^{-w-3L-3}\sum\lambda_{[101],[C,L-2]}^{[-w-3L-2,1]}
((-2\zb{3}\wb{2}W^{C})V)
\\& +(-w-3L-2)(-\mathbf{1})^{-w-3L-3}\sum\lambda_{[110],[C,L-2]}^{[-w-3L-2,1]}
(  
((-2\kappa{}\wb{1}\zb{1}-2\wb{2}\zb{2})W^{C}
)V)
\end{align*}

\paragraph{
Assume \(\kappa+1=0\)} Then \(\SbtSZ{1}{W^C}=0\) and the
expression becomes simpler. 
\begin{align*}
Ldz &=
2 (-w-3L-2)(-\mathbf{1})^{-w-3L-3} (  \sum\lambda_{[010],[C,L]}^{[-w-3L-2,0]}
((-\kappa{}\wb{1}W^{C}
))
 + \sum\lambda_{[100],[C,L]}^{[-w-3L-2,0]}
((\wb{2}W^{C})) )
\\& +2(-w-3L-3)(-\mathbf{1})^{-w-3L-4}\sum\lambda_{[111],[C,L-1]}^{[-w-3L-3,0]}
( 
((\kappa{}\wb{1}\zb{1}\zb{3}+\wb{2}\zb{2}\zb{3})W^{C}
))
%% \\& +(-\mathbf{1})^{-w-3L-2}\sum\lambda_{[000],[C,L-1]}^{[-w-3L-1,1]}
%% ((-w-3L-1)\SbtSZ{\mathbf{1}}{W^{C}}V)
\\& +2 (-w-3L-2)(-\mathbf{1})^{-w-3L-3}\zb{3} (  \sum\lambda_{[011],[C,L-2]}^{[-w-3L-2,1]}
((\kappa{}\wb{1}W^{C}))
 +\sum\lambda_{[101],[C,L-2]}^{[-w-3L-2,1]}
((-\wb{2}W^{C})))V
\\& -2 (-w-3L-2)(-\mathbf{1})^{-w-3L-3}\sum\lambda_{[110],[C,L-2]}^{[-w-3L-2,1]}
(  
((\kappa{}\wb{1}\zb{1}+\wb{2}\zb{2})W^{C}
)V)
\end{align*}
We study the rank of \(Ldz \) by the same way.   
The rank of the first term is 
\( 2 \tbinom{2+L}{L} -  \tbinom{2+L-1}{L-1} \),  
and the third rank is \( 2 \tbinom{2+L-2}{L-2} -  \tbinom{2+L-3}{L-3} \)  
by the same discussion above. 

The second rank is  \(\tbinom{2+L-1}{L-1}\) and   
and the last fourth rank is 
\(  \tbinom{2+L-2}{L-2} \).   
\begin{align*}
\shortintertext{ 
Assume \( -w -m = 2L+1\) and \( 1+\kappa= 0\).  
If  \( -w - 3L -3 < 0\), then }
\dim \pdel \CS{w}{m} = & 0 \\
\shortintertext{If \( -w - 3L -3 = 0\), then }
\dim \pdel \CS{w}{m} = & 
( 2 \tbinom{2+L}{L} -  \tbinom{2+L-1}{L-1} ) +   
( 2 \tbinom{2+L-2}{L-2} -  \tbinom{2+L-3}{L-3} )  
+ 
 \tbinom{2+L-2}{L-2} 
 \\ = &
 2 \tbinom{2+L}{L} -  \tbinom{2+L-1}{L-1}  + 3 
\tbinom{2+L-2}{L-2} 
- \tbinom{2+L-3}{L-3} 
\\
\shortintertext{If \( -w - 3L -3 > 0\), then }
\dim \pdel \CS{w}{m} = & 
( 2 \tbinom{2+L}{L} -  \tbinom{2+L-1}{L-1} ) +   
\tbinom{2+L-1}{L-1} + 
( 2 \tbinom{2+L-2}{L-2} -  \tbinom{2+L-3}{L-3} )  
+ \tbinom{2+L-2}{L-2} 
 \\ = & 
 2 \tbinom{2+L}{L} 
 + 3 \tbinom{2+L-2}{L-2} -  \tbinom{2+L-3}{L-3}   
\end{align*}

% 1+kk ne 0 
\paragraph{
We assume  \( \kappa + 1 \ne 0\)}   

\begin{align*}
Ld &=
(-w-3L-2) (-\mathbf{1})^{-w-3L-3} \zb{3} 
\mw 
\\&\quad 
(\sum\lambda_{[001],[C,L]}^{[-w-3L-2,0]} (\SbtSZ{\mathbf{1}}{W^{C}})
 + \sum\lambda_{[011],[C,L-2]}^{[-w-3L-2,1]} ((2\kappa{}\wb{1}W^{C})V)
+\sum\lambda_{[101],[C,L-2]}^{[-w-3L-2,1]} ((-2\wb{2}W^{C})V)
\\& +(-w-3L-2)(-\mathbf{1})^{-w-3L-3}
\mw \\&\quad 
(  
\sum\lambda_{[010],[C,L]}^{[-w-3L-2,0]}
((-2\kappa{}\wb{1}W^{C}+\zb{2}\SbtSZ{\mathbf{1}}{W^{C}}))
 +\sum\lambda_{[100],[C,L]}^{[-w-3L-2,0]}
((2\wb{2}W^{C}+\zb{1}\SbtSZ{\mathbf{1}}{W^{C}}))
)
\\& +\sum\lambda_{[111],[C,L-1]}^{[-w-3L-3,0]}
( (-\mathbf{1})^{-w-3L-3}          
(2(1+\kappa)\zb{3} V )W^{C}
\\&\qquad 
+(-w-3L-3)(-\mathbf{1})^{-w-3L-4}
((2\kappa{}\wb{1}\zb{1}\zb{3}+2\wb{2}\zb{2}\zb{3})W^{C}+\zb{1}\zb{2}\zb{3}\SbtSZ{\mathbf{1}}{W^{C}}))
%% \\& +(-\mathbf{1})^{-w-3L-2}\sum\lambda_{[000],[C,L-1]}^{[-w-3L-1,1]}
%% ((-w-3L-1)\SbtSZ{\mathbf{1}}{W^{C}}V)
\\& +(-w-3L-2)(-\mathbf{1})^{-w-3L-3}\sum\lambda_{[110],[C,L-2]}^{[-w-3L-2,1]}
(  
((-2\kappa{}\wb{1}\zb{1}-2\wb{2}\zb{2})W^{C}
)V)
\end{align*}
If \( -w -3L-3<0\) then Ld = 0. 
If \( -w -3L-3=0\) then 
\begin{align*}
Ld &=
 \zb{3} \mw 
(\sum\lambda_{[001],[C,L]}^{[-w-3L-2,0]} (\SbtSZ{\mathbf{1}}{W^{C}})
 + \sum\lambda_{[011],[C,L-2]}^{[-w-3L-2,1]} ((2\kappa{}\wb{1}W^{C})V)
+\sum\lambda_{[101],[C,L-2]}^{[-w-3L-2,1]} ((-2\wb{2}W^{C})V)
\\& +
(  
\sum\lambda_{[010],[C,L]}^{[-w-3L-2,0]}
((-2\kappa{}\wb{1}W^{C}+\zb{2}\SbtSZ{\mathbf{1}}{W^{C}}))
 +\sum\lambda_{[100],[C,L]}^{[-w-3L-2,0]}
((2\wb{2}W^{C}+\zb{1}\SbtSZ{\mathbf{1}}{W^{C}}))
)
\\& +\sum\lambda_{[111],[C,L-1]}^{[-w-3L-3,0]}
( (2(1+\kappa)\zb{3} V )W^{C}
)
\\& + \sum\lambda_{[110],[C,L-2]}^{[-w-3L-2,1]}
(  
((-2\kappa{}\wb{1}\zb{1}-2\wb{2}\zb{2})W^{C}
)V)
\\ &=
 \zb{3} \mw 
(\sum\lambda_{[001],[C,L]}^{[-w-3L-2,0]} (\SbtSZ{\mathbf{1}}{W^{C}})
 + \sum\lambda_{[011],[C,L-2]}^{[-w-3L-2,1]} ((2\kappa{}\wb{1}W^{C})V)
 \\& \qquad 
+\sum\lambda_{[101],[C,L-2]}^{[-w-3L-2,1]} ((-2\wb{2}W^{C})V)
+\sum\lambda_{[111],[C,L-1]}^{[-w-3L-3,0]} ( (2(1+\kappa) V )W^{C})
\\& +
(  
\sum\lambda_{[010],[C,L]}^{[-w-3L-2,0]}
((-2\kappa{}\wb{1}W^{C}+\zb{2}\SbtSZ{\mathbf{1}}{W^{C}}))
 +\sum\lambda_{[100],[C,L]}^{[-w-3L-2,0]}
((2\wb{2}W^{C}+\zb{1}\SbtSZ{\mathbf{1}}{W^{C}}))
)
\\& + \sum\lambda_{[110],[C,L-2]}^{[-w-3L-2,1]}
(  
((-2\kappa{}\wb{1}\zb{1}-2\wb{2}\zb{2})W^{C}
)V)
\end{align*}
the first term gives \(\tbinom{2+L-1}{L-1}\), the second term gives 
\( 2\tbinom{2+L}{L} - \tbinom{2+L-1}{L-1} \)  
and the third term gives \( \tbinom{2+L-2}{L-2}\). Thus, the rank is 
\( 2\tbinom{2+L}{L} + \tbinom{2+L-2}{L-2} \).   

If \( -w -3L-3>0\) then the first term gives the rank 
\( \tbinom{2+L-1}{L-1} \), the second term gives 
\( 2\tbinom{2+L}{L} - \tbinom{2+L-1}{L-1} \), the third term gives   
\( \tbinom{2+L-1}{L-1} \), the fourth term gives 
\( \tbinom{2+L-2}{L-2} \). Thus, the rank is 
\[ 
\tbinom{2+L-1}{L-1} + 
 2\tbinom{2+L}{L} - \tbinom{2+L-1}{L-1} + 
\tbinom{2+L-1}{L-1} 
+  \tbinom{2+L-2}{L-2} 
= 2\tbinom{2+L}{L} +  \tbinom{2+L-1}{L-1} +  \tbinom{2+L-2}{L-2} 
\]

%\textcolor{red}{ It's done on Mar 13, 2021.}

\kmcomment{

\CS{w}{-w -2L-1} 
\label{G3:codim:odd}
\\&= 
 \mw ^{-w  -3L -2} \frakg_{-1}
 \mw ^{1} \frakg_{-2}
 \mw ^{ L } \frakg_{-3}
 \mw ^{ 0 } \frakg_{-4}
+
 \mw ^{-w  -3L - 3 } \frakg_{-1}
 \mw ^{ 3 } \frakg_{-2}
 \mw ^{ L  - 1 } \frakg_{-3}
 \mw ^{ 0 } \frakg_{-4}
\notag
\\&\quad 
+
 \mw ^{-w -3L -1 } \frakg_{-1}
 \mw ^{ 0 } \frakg_{-2}
 \mw ^{ L - 1 } \frakg_{-3}
 \mw ^{ 1 } \frakg_{-4}
+
 \mw ^{-w  -3L -2 } \frakg_{-1}
 \mw ^{ 2 } \frakg_{-2}
 \mw ^{ L - 2 } \frakg_{-3}
 \mw ^{ 1 } \frakg_{-4}
\; . \notag
}

\bigskip

\subsubsection{\(\dim \Sbt{\frakg}{\frakg} = 1\) and 
\(\dim \Sbt{\frakg} {\frakg} \not\subset Z(\frakg)\):}

\begin{wrapfigure}[4]{r}[5mm]{0.4\textwidth}
\vspace{-3.0 \baselineskip}
%\begin{center}
%\begin{center}
\( \renewcommand{\arraystretch}{0.8}
\begin{array}{c|c|*{3}{c}|*{3}{c}|c}
\hline
& \mathbf{1} &
\zb{1} &
\zb{2} &
\zb{3} &
\wb{1} &
\wb{2} &
\wb{3} &
V
\\\hline
\mathbf{1} & 0 & 0  & -2\wb{3} & 0 &  -2 V & 0 & 0 &   \\
\cline{1-8}
\zb{1} & 0  & 0 & 0 & 0 & \\
\zb{2} & 2\wb{3} & 0 & 0 &  2 V & \\
\zb{3} & 0 & 0 &  -2 V & 0 &  \\
\cline{1-5}
\wb{1} & -2V &  \\
\wb{2} & 0 & \\
\wb{3} & 0 &  \\
\cline{1-2}
V   \\
\end{array}
\)
%\end{center}
\end{wrapfigure}
In the case of \(\dim \Sbt{\frakg}{\frakg} = 1\) and \(\dim \Sbt{\frakg}
{\frakg} \not\subset Z(\frakg)\), we may take a basis \( \Sbt{\xi_{1}}
{\xi_{2}} = - \Sbt{\xi_{2}}{\xi_{1}} = \xi_{2}\) and the others are 0.  Now
we have \( d \zb{1} =  0 \;,\;  d \zb{2} =  -2 \ww{3}\;,\; d \zb{3} = 0 \).
The brackets of superalgebra are given as follows:

\eqref{bdary:d3} implies 
\begin{align*}
\pdel ( \mathbf{1}^a Z^B W^C V^{\ell}) &= 
 (-\mathbf{1})^a (\pdel Z^{B}) W^C V^{\ell} 
 + a (-\mathbf{1})^{a-1}( \SbtSZ{ \mathbf{1} } {Z^{B}} W^C 
 - \SbtSZ{ \mathbf{1} }{ W^C } Z^{B} ) V^{\ell} 
\\\shortintertext{where}
\pdel Z^{111} & = -2 \zb{1}\mw V\;,\quad  
 \pdel Z^{110} =  \pdel Z^{101} = 0\;, \quad    
  \pdel Z^{011} = 2 V \;,    
\\
\SbtSZ{ \mathbf{1} } {Z^{B}} &= 
 2\parity{1+b_{3}} \frac{ b_{2}}{\zb{2}} Z^{B} \wb{3} 
\\
\SbtSZ{ \mathbf{1} }{ W^C } &=
   ( -2 V )  \frac{ c_{1} } {\wb{1}} W^{C}   
\\\shortintertext{ 
because}
\SbtSZ{ \mathbf{1} } {Z^{B}} &= 
b_{1} \SbtSZ{\mathbf{1}}{\zb{1}} \zb{2}^{b_{2}}\zb{3}^{b_{3}}
+ b_{2} \zb{1}^{b_{1}} \SbtSZ{\mathbf{1}}{\zb{2}}\zb{3}^{b_{3}}
+ b_{3} \zb{1}^{b_{1}} \zb{2}^{b_{2}} \SbtSZ{\mathbf{1}}{\zb{3}}
= b_{2} \zb{1}^{b_{1}}( -2 \wb{3} ) \zb{3}^{b_{3}}
\\& 
= 2 b_{2} \parity{1+b_{3}} \zb{1}^{b_{1}} \zb{3}^{b_{3}} \wb{3} 
= 2\parity{1+b_{3}} \frac{ b_{2}}{\zb{2}} Z^{B} \wb{3} 
\\
\SbtSZ{ \mathbf{1} }{ W^C } &=
 c_{1} \SbtSZ{ \mathbf{1} }{\wb{1} }  
\wb{1}^{c_{1}-1} \wb{2}^{c_{2}} \wb{3}^{c_{3}} 
=  
  \SbtSZ{ \mathbf{1} }{\wb{1} } \frac{ c_{1} } {\wb{1}} W^{C}  
=   ( -2 V )  \frac{ c_{1} } {\wb{1}} W^{C}   
\kmcomment{
\\
\SbtSZ{ \mathbf{1}^{a} }{ \zb{i} W^C } &=
a (- \mathbf{1})^{a-1} ( -2 \delta_{i}^{2} \wb{3} W^{C}
 + 2 V  \frac{ c_{1} } {\wb{1}} W^{C})  
 \\
\SbtSZ{ \mathbf{1}^{a} }{ \zb{i}\zb{j} W^C } &=
a (- \mathbf{1})^{a-1} ( -2 \begin{vmatrix} \delta^{2}_{i} &  \delta^{2}_{j}
\\ \zb{i}& \zb{j} \end{vmatrix} \wb{3} W^{C} + 2 \zb{i}\zb{j} V  \frac{ c_{1} } {\wb{1}} W^{C})   
\\
\SbtSZ{ \mathbf{1}^{a} }{ \zb{1}\zb{2}\zb{3} W^C } &=
a (- \mathbf{1})^{a-1} ( 2  \zb{1} \zb{3} \wb{3} W^{C} + 2 \parity{|C|} 
\zb{1}\zb{2}\zb{3} V  \frac{ c_{1} } {\wb{1}} W^{C})   
}
\end{align*}

Assume \(-w -m =2K\). Then a general element in  \(\pdel \CS{w}{m}\) is
given by 
\begin{align*}
%ukeAllKc
X &=\sum\lambda_{[000],[C,K]}^{[-w-3K,0]}((-w-3K)(-\mathbf{1})^{-w-3K-1}\SbtSZ{\mathbf{1}}{W^{C}})
\\&
+\sum\lambda_{[011],[C,K-1]}^{[-w-3K-1,0]}(2(-\mathbf{1})^{-w-3K-1} V W^{C}
\\&\quad 
+(-w-3K-1)(-\mathbf{1})^{-w-3K-2}(2\zb{3}\wb{3}W^{C}+\zb{2}\zb{3}\SbtSZ{\mathbf{1}}{W^{C}}))
\\&
+\sum\lambda_{[101],[C,K-1]}^{[-w-3K-1,0]}((-w-3K-1)(-\mathbf{1})^{-w-3K-2}\zb{1}\zb{3}\SbtSZ{\mathbf{1}}{W^{C}})
\\&
+\sum\lambda_{[110],[C,K-1]}^{[-w-3K-1,0]}((-w-3K-1)(-\mathbf{1})^{-w-3K-2}(-2\zb{1}\wb{3}W^{C}+\zb{1}\zb{2}\SbtSZ{\mathbf{1}}{W^{C}}))
\\&
+\sum\lambda_{[001],[C,K-2]}^{[-w-3K,1]}( 0 )
\\&
+\sum\lambda_{[010],[C,K-2]}^{[-w-3K,1]}((-w-3K)(-\mathbf{1})^{-w-3K-1}
(-2\wb{3}W^{C}+ 0 )V)
\\&
+\sum\lambda_{[100],[C,K-2]}^{[-w-3K,1]}( 0 ) 
\\&
+\sum\lambda_{[111],[C,K-3]}^{[-w-3K-1,1]}( 0
+(-w-3K-1)(-\mathbf{1})^{-w-3K-2}(2\zb{1}\zb{3}\wb{3}W^{C}+ 0)V)
\\ =& 
+ (-w-3K)(-\mathbf{1})^{-w-3K-1}
\sum\lambda_{[000],[C,K]}^{[-w-3K,0]}(\SbtSZ{\mathbf{1}}{W^{C}})
\\&
+\sum\lambda_{[011],[C,K-1]}^{[-w-3K-1,0]}(2(-\mathbf{1})^{-w-3K-1} V W^{C}
\\&\quad 
+(-w-3K-1)(-\mathbf{1})^{-w-3K-2}(2\zb{3}\wb{3}W^{C}+\zb{2}\zb{3}\SbtSZ{\mathbf{1}}{W^{C}}))
\\&
+ (-w-3K-1)(-\mathbf{1})^{-w-3K-2}
\sum\lambda_{[101],[C,K-1]}^{[-w-3K-1,0]}(\zb{1}\zb{3}\SbtSZ{\mathbf{1}}{W^{C}})
\\&
+\sum\lambda_{[110],[C,K-1]}^{[-w-3K-1,0]}((-w-3K-1)(-\mathbf{1})^{-w-3K-2}(-2\zb{1}\wb{3}W^{C}+\zb{1}\zb{2}\SbtSZ{\mathbf{1}}{W^{C}}))
\\&
+(-w-3K)(-\mathbf{1})^{-w-3K-1}
\sum\lambda_{[010],[C,K-2]}^{[-w-3K,1]}(
(-2\wb{3}W^{C})V)
\\&
+ (-w-3K-1)(-\mathbf{1})^{-w-3K-2}
\sum\lambda_{[111],[C,K-3]}^{[-w-3K-1,1]}((2\zb{1}\zb{3}\wb{3}W^{C})V)
\\ =& 
+ (-w-3K)(-\mathbf{1})^{-w-3K-1} (
\sum\lambda_{[000],[C,K]}^{[-w-3K,0]}\SbtSZ{\mathbf{1}}{W^{C}}
+ \sum\lambda_{[010],[C,K-2]}^{[-w-3K,1]} (-2\wb{3}W^{C})V )
\\& 
+\sum\lambda_{[011],[C,K-1]}^{[-w-3K-1,0]}(2(-\mathbf{1})^{-w-3K-1} V W^{C}
\\&\quad 
+(-w-3K-1)(-\mathbf{1})^{-w-3K-2}(2\zb{3}\wb{3}W^{C}+\zb{2}\zb{3}\SbtSZ{\mathbf{1}}{W^{C}}))
\\& 
+ (-w-3K-1)(-\mathbf{1})^{-w-3K-2} \zb{1}\zb{3}(
\sum\lambda_{[101],[C,K-1]}^{[-w-3K-1,0]}\SbtSZ{\mathbf{1}}{W^{C}}
%\\& 6
+ \sum\lambda_{[111],[C,K-3]}^{[-w-3K-1,1]}2\wb{3}W^{C}V)
\\& 
+(-w-3K-1)(-\mathbf{1})^{-w-3K-2}
\sum\lambda_{[110],[C,K-1]}^{[-w-3K-1,0]}((-2\zb{1}\wb{3}W^{C}+\zb{1}\zb{2}\SbtSZ{\mathbf{1}}{W^{C}}))
\end{align*}
If \( -w -3K-1 < 0\), then \(X=0\), i.e., the boundary operator is zero map. 

If \( -w -3K-1 = 0\),  
\begin{align*} X =& 
 ( \sum\lambda_{[000],[C,K]}^{[-w-3K,0]}\SbtSZ{\mathbf{1}}{W^{C}}
+ \sum\lambda_{[010],[C,K-2]}^{[-w-3K,1]} (-2\wb{3}W^{C})V )
+\sum\lambda_{[011],[C,K-1]}^{[-w-3K-1,0]} 2 V W^{C} 
\\=& 
2V (- \sum\lambda_{[000],[C,K]}^{[-w-3K,0]}
     \frac{ c_{1} } {\wb{1}} W^{C}   
% \SbtSZ{\mathbf{1}}{W^{C}}
+\parity{K} \sum\lambda_{[010],[C,K-2]}^{[-w-3K,1]}  \wb{3}W^{C} 
+\sum\lambda_{[011],[C,K-1]}^{[-w-3K-1,0]} W^{C}) 
\\=& 
2V \sum_{|C|=K-1} (- \lambda_{[000],[[1+c_1,c_2,c_3],K]}^{[-w-3K,0]}
     (1+ c_{1} ) 
+\parity{K} \lambda_{[010],[[c_1,c_2, -1+c_3],K-2]}^{[-w-3K,1]}  
+\lambda_{[011],[[c_1,c_2,c_3],K-1]}^{[-w-3K-1,0]} )W^{C} 
\end{align*}
and we have the set of generators:
\[
 - \lambda_{[000],[[1+c_1,c_2,c_3],K]}^{[-w-3K,0]}
     (1+ c_{1} ) 
+\parity{K} \lambda_{[010],[[c_1,c_2, -1+c_3],K-2]}^{[-w-3K,1]}  
+\lambda_{[011],[[c_1,c_2,c_3],K-1]}^{[-w-3K-1,0]} 
\] with \(|C|=K-1\), and the rank is  \( \tbinom{2+K-1}{K-1}\).   

If \( -w -3K-1 > 0\), we have linearly independent 4 blocks:
\begin{align*}
X_ 1 =& 
+ (-w-3K)(-\mathbf{1})^{-w-3K-1} (
\sum\lambda_{[000],[C,K]}^{[-w-3K,0]}\SbtSZ{\mathbf{1}}{W^{C}}
+ \sum\lambda_{[010],[C,K-2]}^{[-w-3K,1]} (-2\wb{3}W^{C})V )
\\
X_ 2 = & 
+\sum\lambda_{[011],[C,K-1]}^{[-w-3K-1,0]}(2(-\mathbf{1})^{-w-3K-1} V W^{C}
\\&\quad 
+(-w-3K-1)(-\mathbf{1})^{-w-3K-2}(2\zb{3}\wb{3}W^{C}+\zb{2}\zb{3}\SbtSZ{\mathbf{1}}{W^{C}}))
\\
X_3 = & 
+ (-w-3K-1)(-\mathbf{1})^{-w-3K-2} \zb{1}\zb{3}(
\sum\lambda_{[101],[C,K-1]}^{[-w-3K-1,0]}\SbtSZ{\mathbf{1}}{W^{C}}
%\\& 6
+ \sum\lambda_{[111],[C,K-3]}^{[-w-3K-1,1]}2\wb{3}W^{C}V)
\\
X_4 = & 
+\sum\lambda_{[110],[C,K-1]}^{[-w-3K-1,0]}((-w-3K-1)(-\mathbf{1})^{-w-3K-2}
(-2\zb{1}\wb{3}W^{C}+\zb{1}\zb{2}\SbtSZ{\mathbf{1}}{W^{C}})) \;. 
\end{align*}
Only 
\(X_2\) and \( X_4\) have terms which are free from \(V\) and linearly
independent each other. Their freedom is \( 2 \tbinom{2+K-1}{K-1}\).   
% About  \(X_1\) and \( X_3\), we have generators
Since 
\begin{align*}
 X_ 1 =& 
(-w-3K)(-\mathbf{1})^{-w-3K-1} ( -2V) (  
 \sum\lambda_{[000],[C,K]}^{[-w-3K,0]}
     \frac{ c_{1} } {\wb{1}} W^{C}   
- \parity{ K } \sum\lambda_{[010],[C,K-2]}^{[-w-3K,1]}\wb{3}W^{C} 
)
\\
X_3 = & 
 (-w-3K-1)(-\mathbf{1})^{-w-3K-2} \zb{1}\zb{3}( -2 V) (
\sum\lambda_{[101],[C,K-1]}^{[-w-3K-1,0]} \frac{ c_{1} } {\wb{1}} W^{C}   
%\\& 6
-\parity{K} \sum\lambda_{[111],[C,K-3]}^{[-w-3K-1,1]}\wb{3}W^{C})
\\\shortintertext{
we see that generators without common factors are }
& 
 \lambda_{[000],[[1+c_1,c_2,c_3],K]}^{[-w-3K,0]}
     ( 1+ c_{1} ) 
- \parity{ K } \sum\lambda_{[010],[[c_1,c_2,-1+c_3],K-2]}^{[-w-3K,1]}
\quad \text{ with } c_1+c_2+c_3 = K-1
\\ 
& 
\lambda_{[101],[[1+c_1,c_2,c_3],K-1]}^{[-w-3K-1,0]} (1+ c_{1} ) 
%\\& 6
-\parity{K} \lambda_{[111],[[c_1,c_2, -1+c_3],K-3]}^{[-w-3K-1,1]}
\quad \text{ with } c_1+c_2+c_3 = K-2 \;.
\end{align*}
Their freedom are 
 \(  \tbinom{2+K-1}{K-1}\) and  \(  \tbinom{2+K-2}{K-2}\) separately. Thus, the rank is    
 \( 3 \tbinom{2+K-1}{K-1} +   \tbinom{2+K-2}{K-2}\).  

\medskip

When 
\(-w -m =2L+1\),  a general element in \(\pdel \CS{w}{m}\) is given by 
\begin{align*}
Y =& + \sum\lambda_{[001],[C,L]}^{[-w-3L-2,0]}((-w-3L-2)(-\mathbf{1})^{-w-3L-3}\zb{3}\SbtSZ{\mathbf{1}}{W^{C}})
\\&
+\sum\lambda_{[010],[C,L]}^{[-w-3L-2,0]}((-w-3L-2)(-\mathbf{1})^{-w-3L-3}(-2\wb{3}W^{C}+\zb{2}\SbtSZ{\mathbf{1}}{W^{C}}))
\\&
+\sum\lambda_{[100],[C,L]}^{[-w-3L-2,0]}((-w-3L-2)(-\mathbf{1})^{-w-3L-3}\zb{1}\SbtSZ{\mathbf{1}}{W^{C}})
\\&
+\sum\lambda_{[111],[C,L-1]}^{[-w-3L-3,0]}(-2(-\mathbf{1})^{-w-3L-3}\zb{1} V W^{C}
\\&\quad 
+(-w-3L-3)(-\mathbf{1})^{-w-3L-4}(2\zb{1}\zb{3}\wb{3}W^{C}+\zb{1}\zb{2}\zb{3}\SbtSZ{\mathbf{1}}{W^{C}}))
\\&
+\sum\lambda_{[000],[C,L-1]}^{[-w-3L-1,1]}( 0 )
\\&
+\sum\lambda_{[011],[C,L-2]}^{[-w-3L-2,1]}( 0 
+(-w-3L-2)(-\mathbf{1})^{-w-3L-3}(2\zb{3}\wb{3}W^{C}+ 0)V)
\\&
+\sum\lambda_{[101],[C,L-2]}^{[-w-3L-2,1]}( 0 )
\\&
+\sum\lambda_{[110],[C,L-2]}^{[-w-3L-2,1]}((-w-3L-2)(-\mathbf{1})^{-w-3L-3}
(-2\zb{1}\wb{3}W^{C}+ 0 )V)\; . 
\\
=& + (-w-3L-2)(-\mathbf{1})^{-w-3L-3}\zb{3} (
\sum\lambda_{[001],[C,L]}^{[-w-3L-2,0]}(\SbtSZ{\mathbf{1}}{W^{C}})
+
\sum\lambda_{[011],[C,L-2]}^{[-w-3L-2,1]}(2\wb{3}W^{C}V) )
\\& 
+(-w-3L-2)(-\mathbf{1})^{-w-3L-3}
\sum\lambda_{[010],[C,L]}^{[-w-3L-2,0]}(-2\wb{3}W^{C}+\zb{2}\SbtSZ{\mathbf{1}}{W^{C}})
\\&
+
\sum\lambda_{[111],[C,L-1]}^{[-w-3L-3,0]}(-2(-\mathbf{1})^{-w-3L-3}\zb{1} V W^{C}
\\&\quad 
+(-w-3L-3)(-\mathbf{1})^{-w-3L-4}(2\zb{1}\zb{3}\wb{3}W^{C}+\zb{1}\zb{2}\zb{3}\SbtSZ{\mathbf{1}}{W^{C}}))
\\&
+(-w-3L-2)(-\mathbf{1})^{-w-3L-3}\zb{1}
( \sum\lambda_{[100],[C,L]}^{[-w-3L-2,0]}\SbtSZ{\mathbf{1}}{W^{C}}
-
\sum\lambda_{[110],[C,L-2]}^{[-w-3L-2,1]} 2\wb{3}W^{C} V)\; . 
\end{align*}
By the same discussion when \(-w-m = 2K\), we also get \( \dim \CS{w}{m}\)
for \(-w-m = 2L+1\) as follows: 

\begin{prop}
When \( -w -m=2K\), then 
\( \dim \pdel \CS{w}{m} = \begin{cases}
  3 \tbinom{2+K-1}{K-1} +   \tbinom{2+K-2}{K-2} & \text{if } -w -3K-1 > 0 \\
 \tbinom{2+K-1}{K-1}  &   \text{if }-w -3K-1 = 0  \\  
0 & \text{otherwise.}
\end{cases}\)

When \( -w -m=2L+1\), then 
\( \dim \pdel \CS{w}{m} = \begin{cases}
  3 \tbinom{2+L-1}{L-1} +   \tbinom{2+L}{L} & \text{if } -w -3L-2 > 1 \\
  2 \tbinom{2+L-1}{L-1} +   \tbinom{2+L}{L} & \text{if } -w -3L-2 = 1 \\
0 &\text{otherwise.}  
\end{cases}\)
\end{prop}

\bigskip

%%%%%%%%%%%%%%%%%%%%%%%%%%%%%%%%%%%%%%%%%%%%%%%%%%%%%%%%%%%%%%%%%
%%  \(\dim \Sbt{\frakg}{\frakg} = 1\) and 
%%  \(\Sbt{\frakg}{\frakg} \subset Z(\frakg)\), 
%%%%%%%%%%%%%%%%%%%%%%%%%%%%%%%%%%%%%%%%%%%%%%%%%%%%%%%%%%%%%%%%%
\subsubsection{\(\dim \Sbt{\frakg}{\frakg} = 1\) and 
 \(\Sbt{\frakg}{\frakg} \subset Z(\frakg)\):}
\kmcomment{
In the case of \(\dim \Sbt{\frakg}{\frakg} = 1\) and 
 \(\Sbt{\frakg}{\frakg} \subset Z(\frakg)\), we may take a basis 
\( \Sbt{\xi_{1}}{\xi_{2}} =
- \Sbt{\xi_{2}}{\xi_{1}} = \xi_{3}\) and the others are 0. 
Now we have \( d \zb{1} =  0 \;,\;  d \zb{2} =  0 \;,\;
d \zb{3} = -2 \wb{3} \).  The brackets of superalgebra are given as follows:
}
\begin{wrapfigure}[4]{r}[5mm]{0.4\textwidth}
\vspace{-2.8 \baselineskip}
%\begin{center}
\( \renewcommand{\arraystretch}{0.8}
\begin{array}{c|c|*{3}{c}|*{3}{c}|c}
\hline
& \mathbf{1} &
\zb{1} &
\zb{2} &
\zb{3} &
\wb{1} &
\wb{2} &
\wb{3} &
V
\\\hline
\mathbf{1} & 0 & 0  & 0 &  -2\wb{3}  &  0 & 0 & 0 &   \\
\cline{1-8}
\zb{1} & 0  & 0 & 0 & 0 & \\
\zb{2} & 0 & 0 & 0 &   0 & \\
\zb{3} & 2\wb{3}  & 0 &  0  & 0 &  \\
\cline{1-5}
\wb{1} & 0 &  \\
\wb{2} & 0 & \\
\wb{3} & 0 &  \\
\cline{1-2}
V   \\
\end{array}
\)
%\end{center}
\end{wrapfigure}
In the case of \(\dim \Sbt{\frakg}{\frakg} = 1\) and 
 \(\Sbt{\frakg}{\frakg} \subset Z(\frakg)\), we may take a basis 
\( \Sbt{\xi_{1}}{\xi_{2}} =
- \Sbt{\xi_{2}}{\xi_{1}} = \xi_{3}\) and the others are 0. 
Now we have \( d \zb{1} =  0 \;,\;  d \zb{2} =  0 \;,\;
d \zb{3} = -2 \wb{3} \).  The brackets of superalgebra are given as on the
right

\eqref{bdary:d3} implies 
\begin{align*}
\pdel ( \mathbf{1}^a Z^B W^C V^{\ell}) &= 
 a (- \mathbf{1})^{a-1} \SbtSZ{ \mathbf{1} } {Z^{B}} W^C V^{\ell} 
= a (- \mathbf{1})^{a-1} \zb{1}^{b_{1}} \zb{2}^{b_{2}} b_{3}(-2 \wb{3})
 W^C V^{\ell}
 \\& 
= a (- \mathbf{1})^{a-1}  b_{3} \frac{ Z^{B} }{\zb{3}}(-2 \wb{3})
 W^C V^{\ell}
 \;.  
\end{align*}
\begin{align*}
X &=\sum\lambda_{[000],[C,K]}^{[-w-3K,0]}(0)
+\sum\lambda_{[011],[C,K-1]}^{[-w-3K-1,0]}(-2(-w-3K-1)(-\mathbf{1})^{-w-3K-2}\zb{2} \wb{3} W^{C})
\\&
+\sum\lambda_{[101],[C,K-1]}^{[-w-3K-1,0]}(-2(-w-3K-1)(-\mathbf{1})^{-w-3K-2}\zb{1} \wb{3} W^{C})
+\sum\lambda_{[110],[C,K-1]}^{[-w-3K-1,0]}(0)
\\&
+\sum\lambda_{[001],[C,K-2]}^{[-w-3K,1]}(-2(-w-3K)(-\mathbf{1})^{-w-3K-1} \wb{3} W^{C}V)
\\&
+\sum\lambda_{[010],[C,K-2]}^{[-w-3K,1]}(0)
+\sum\lambda_{[100],[C,K-2]}^{[-w-3K,1]}(0)
\\&
+\sum\lambda_{[111],[C,K-3]}^{[-w-3K-1,1]}(-2(-w-3K-1)(-\mathbf{1})^{-w-3K-2}\zb{1}\zb{2} \wb{3} W^{C}V)
\\
Y &=\sum\lambda_{[001],[C,L]}^{[-w-3L-2,0]}(-2(-w-3L-2)(-\mathbf{1})^{-w-3L-3} \wb{3} W^{C})
\\& 
+\sum\lambda_{[010],[C,L]}^{[-w-3L-2,0]}(0)
+\sum\lambda_{[100],[C,L]}^{[-w-3L-2,0]}(0)
\\&
+\sum\lambda_{[111],[C,L-1]}^{[-w-3L-3,0]}(-2(-w-3L-3)(-\mathbf{1})^{-w-3L-4}\zb{1}\zb{2} \wb{3} W^{C})
+\sum\lambda_{[000],[C,L-1]}^{[-w-3L-1,1]}(0)
\\&
+\sum\lambda_{[011],[C,L-2]}^{[-w-3L-2,1]}(-2(-w-3L-2)(-\mathbf{1})^{-w-3L-3}\zb{2} \wb{3} W^{C}V)
\\&
+\sum\lambda_{[101],[C,L-2]}^{[-w-3L-2,1]}(-2(-w-3L-2)(-\mathbf{1})^{-w-3L-3}\zb{1} \wb{3} W^{C}V)
+\sum\lambda_{[110],[C,L-2]}^{[-w-3L-2,1]}(0)\;. 
\end{align*}
Thus, when \(-w -m =2K\), \(\pdel \CS{w}{m}\) is spanned by   
\begin{align*}
& 
(-\mathbf{1})^{-w -3K-2} \zb{2}( -2 \wb{3} ) W^{C} \text{ and } 
 (-\mathbf{1})^{-w -3K-2} \zb{1}( -2 \wb{3} ) W^{C}  \text{ with } |C|=K-1 \; ,
\\& 
(-\mathbf{1})^{-w -3K-1} ( -2 \wb{3}) W^{C}V  \text{ with } |C|=K-2\;, 
\\& 
(-\mathbf{1})^{-w -3K-2} \zb{1}\zb{2}(-2 \wb{3}) W^{C}V  \text{ with } |C|=K-3 
\end{align*}
and \(\dim \pdel \CS{w}{m}
= \tbinom{-w-3K-2}{0} 2 \tbinom{2+K-1}{K-1}
+ \tbinom{-w-3K-1}{0} \tbinom{2+K-2}{K-2}
+ \tbinom{-w-3K-2}{0} \tbinom{2+K-3}{K-3}
\).  

By the same argument, 
if \(-w -m =2L+1\), then 
\[\dim \pdel \CS{w}{m} =  
\tbinom{-w-3L-3}{0}  \tbinom{2+L}{L}
+ \tbinom{-w-3L-4}{0}  \tbinom{2+L-1}{L-1}
+ \tbinom{-w-3L-3}{0} 2 \tbinom{2+L-2}{L-2}
\;.
\]
The $m$-th 
Betti number is given by 
\( \dim \CS{w}{m} - \dim \pdel \CS{w}{m} - \dim \pdel \CS{w}{m+1} \), say
\(\text{Bet}_{m} \). Thus,    
\begin{align*}
\shortintertext{
if \( -w -m = 2K\) then} 
\text{Bet}_{m} & = \dim \CS{w}{m} - 
( \tbinom{-w-3K-2}{0} 2 \tbinom{2+K-1}{K-1}
+ \tbinom{-w-3K-1}{0} \tbinom{2+K-2}{K-2}
+ \tbinom{-w-3K-2}{0} \tbinom{2+K-3}{K-3}
)
\\& 
- 
(
\tbinom{-w-3L}{0}  \tbinom{2+L+1}{L-1}
+ \tbinom{-w-3L-1}{0}  \tbinom{2+L}{L-2}
+ \tbinom{-w-3L}{0} 2 \tbinom{2+L-1}{L-3}
)
\\
\shortintertext{
and 
if \( -w -m = 2L+1 \) then} 
\text{Bet}_{m}& = \dim \CS{w}{m} 
- 
(
\tbinom{-w-3L-3}{0}  \tbinom{2+L}{L}
+ \tbinom{-w-3L-4}{0}  \tbinom{2+L-1}{L-1}
+ \tbinom{-w-3L-3}{0} 2 \tbinom{2+L-2}{L-2}
)
\\& 
-
(
 \tbinom{-w-3L-2}{0} 2 \tbinom{2+L-1}{L-1}
+ \tbinom{-w-3L-1}{0} \tbinom{2+L-2}{L-2}
+ \tbinom{-w-3L-2}{0} \tbinom{2+L-3}{L-3}
)
\end{align*}

\kmcomment{
1 & 1\in \frakg_{-1} & \sigma_{i} \in \frakg_{-2} 
 & \sigma_{j}\wedge \sigma_{k} \in \frakg_{-3} 
 &  \sigma_{1}\wedge   \sigma_{2}\wedge   \sigma_{3}
 \in \frakg_{-4} 
 \\\hline
}

\subsubsection{Look at once}
So far, for each chain complex of Lie  superalgebra of 
a given  3-dimensional Lie algebra, we got all kernel dimensions, and so 
we can get Betti numbers.  For instance, the following table shows 
the case of weight \(-3, -5, -10\) chain spaces simultaneously. 
\(d3:ker\) means the list of kernel dimensions of the Lie algebra
\(\dim \Sbt{\frakg}{\frakg} = 3\), 
\(d2y:\) means the Lie algebra
\(\dim \Sbt{\frakg}{\frakg} = 2\) and \( \kappa +1 = 0\), 
\(d2n:\) means the Lie algebra
\(\dim \Sbt{\frakg}{\frakg} = 2\) and \( \kappa +1 \ne 0\), 
\(d1n:\) means the Lie algebra \(\dim \Sbt{\frakg}{\frakg} = 1\) and 
\(\Sbt{\frakg}{\frakg} \not\subset Z( \frakg )\),  and 
\(d1y:\) means the Lie algebra \(\dim \Sbt{\frakg}{\frakg} = 1\) and 
\(\Sbt{\frakg}{\frakg} \subset Z( \frakg )\).

\begin{center}
\(
\begin{array}{c| *{3}{c}}
w = -3 & 1 & 2 & 3 \\\hline
SpDim & 
3 & 3 & 1 \\\hline
d3:\ker &  
3 & 0 & 1
\\\phantom{d3:}Bet  &  
0 & 0 & 1 \\\hline
d2y:\ker & 
3 & 1 & 1
\\\phantom{d2y:}Bet  & 
1 & 1 & 1 \\\hline
d2n:\ker & 
3 & 1 & 1
\\\phantom{d2n:}Bet  & 
1 & 1 & 1 \\\hline
d1y:\ker &
3 & 2 & 1
\\\phantom{ d1y:}Bet  &
2 & 2 & 1 \\\hline
d1n:\ker & 
3 & 2 & 1
\\\phantom{d1n:}Bet  & 
2 & 2 & 1 \\\hline
\end{array}
\)
\hfil 
\(
\begin{array}{c|*{5}{c}}
w=-5 & 1 & 2 & 3 & 4 & 5 \\\hline
SpDim & 
0 & 10 & 6 & 3 & 1 \\\hline
d3 :\ker &
0 & 10 & 3 & 0 & 1 \\\phantom{d3:}Bet & 
0 & 7 & 0 & 0 & 1 \\\hline
d2y :\ker &
0 & 10 & 3 & 1 & 1 \\\phantom{d2y:}Bet & 
0 & 7 & 1 & 1 & 1 \\\hline
d2n :\ker &
0 & 10 & 2 & 1 & 1 \\\phantom{d2n:}Bet & 
0 & 6 & 0 & 1 & 1 \\\hline
d1y :\ker &
0 & 10 & 4 & 2 & 1 \\\phantom{d1y:}Bet & 
0 & 8 & 3 & 2 & 1 \\\hline
d1n :\ker &
0 & 10 & 3 & 2 & 1 \\\phantom{d1n:}Bet & 
0 & 7 & 2 & 2 & 1 \\\hline
\end{array}
\)
\hfil
\(
\begin{array}{c| *{10}{c}}
\kmcomment{Bango}
w = - 10 & 1 & 2 & 3 & 4 & 5 & 6 & 7 & 8 & 9 & 10 \\\hline
Space\; \dim &
0 & 0 & 6 & 38 & 27 & 18 & 11 & 6 & 3 & 1 \\\hline
d3:\ker &
0 & 0 & 6 & 32 & 11 & 7 & 4 & 3 & 0 & 1 \\
\phantom{d3:}Bet & 
0 & 0 & 0 & 16 & 0 & 0 & 1 & 0 & 0 & 1 \\\hline
d2y:\ker &
0 & 0 & 6 & 33 & 12 & 8 & 5 & 3 & 1 & 1 \\
\phantom{d2y:}Bet &
0 & 0 & 1 & 18 & 2 & 2 & 2 & 1 & 1 & 1 \\\hline
d2n:\ker &
0 & 0 & 6 & 32 & 11 & 7 & 4 & 2 & 1 & 1 \\
\phantom{d2n:}Bet &
0 & 0 & 0 & 16 & 0 & 0 & 0 & 0 & 1 & 1 \\\hline
d1y:\ker &
0 & 0 & 6 & 35 & 16 & 11 & 7 & 4 & 2 & 1 \\
\phantom{d1y:}Bet &
0 & 0 & 3 & 24 & 9 & 7 & 5 & 3 & 2 & 1 \\\hline
d1n:\ker & 
0 & 0 & 6 & 32 & 12 & 8 & 5 & 3 & 2 & 1 \\
\phantom{d1n:}Bet & 
0 & 0 & 0 & 17 & 2 & 2 & 2 & 2 & 2 & 1 \\\hline
\end{array}
\)
\end{center}

\section{Long \(\mZ\)-graded Lie superalgebra }

\subsection{Examples of \(\mZ\)-graded trivially long Lie superalgebra}
\kmcomment{
\(\oplus_{i} \Lambda^{i} \cbdl{M}\oplus \oplus_{j} \Lambda^{j} \tbdl{M}\) 
}
Let \(M\) be an \(n\)-dimensional manifold. 
It is known that 
\begin{align}
&  \Lambda ^{ 0 } \tbdl{M}\oplus 
 \Lambda ^{ 1 } \tbdl{M}\oplus \cdots  \oplus 
 \Lambda ^{ n } \tbdl{M}\qquad \text{ grade of }\Lambda^{i}\tbdl{M} = i-1
\label{super:tbdl}
 \\
\shortintertext{is a Lie superalgebra with the Schouten bracket with grading  
of \(\Lambda^{i}\tbdl{M} = i-1\) and  \( \SbtS{f}{g}= 0\) for functions.   
In this article, we have shown that }
&
 \Lambda ^{ n } \cbdl{M}\oplus \cdots  \oplus 
 \Lambda ^{ 1 } \cbdl{M}\oplus 
 \Lambda ^{ 0 } \cbdl{M}\qquad \text{ grade of }\Lambda^{j}\cbdl{M} = -1-j
\label{super:cbdl}
 \\
 \shortintertext{ 
becomes a Lie superalgebra  with the bracket defined
by \eqref{defn:bkt:forms}.  
Now let} 
 \frakg & = 
 \Lambda ^{ n } \cbdl{M}\oplus \cdots\oplus   
 \Lambda ^{ 1 } \cbdl{M}\oplus 
 \Lambda ^{ 0 } \cbdl{M}\oplus 
 \Lambda ^{ 1 } \tbdl{M}\oplus  \cdots\oplus   
 \Lambda ^{ n } \tbdl{M} \; . 
\label{super:mixed}
 \end{align} 
We easily get a ``long trivial'' \(\mZ\)-graded Lie superalgebra as follows:

\begin{prop} 
Define the grading of \eqref{super:mixed} by  
\( \frakg_{i} = \Lambda ^{ -1-i } \cbdl{M}\) for \( i <  0\), and 
 \( \frakg_{i} = \Lambda ^{ 1+i } \tbdl{M}\)  for \( i \geqq 0\), i.e., 
 \( \Lambda ^{ j } \tbdl{M} = \frakg_{-1 + j}\) for \( j>0\),  and \(  
 \Lambda ^{ j } \cbdl{M} = 
\frakg_{-(1+j)} \) for \( j\geqq 0\).  

\kmcomment{
\( \frakg_{i} = \begin{cases} \Lambda ^{ 1+i } \tbdl{M} & \text{for } i\geqq 0 \\
 \Lambda ^{ -1-i } \cbdl{M} &\text{otherwise} \end{cases}\), 
}

On \( \frakg \), define a bracket by  
\( \SbtE{new}{x}{y} = \begin{cases} \SbtS{x}{y} & x,y \text{ vectors} \\
\SbtD{x}{y} & x,y \text{ forms} \\ 0 & \text{otherwise} \end{cases}\). 
Then \(\frakg\) becomes a trivial 
Lie superalgebra. 
\end{prop}
\kmcomment{
Comparing the previous Proposition, here comes a simple and natural question: 
\begin{quote}
Is there any non-trivial bracket for
forms and multivector, by which, the direct sum \(\ds \mathop{\oplus}_{i}
\Lambda^{i} \cbdl{M}  \mathop{\oplus}_{j} \Lambda^{j} \tbdl{M}\)
becomes Lie superalgebra?  
\end{quote}
}
We restate the proposition above in 
 abstract way as below:  
\begin{prop}
Let \( \frakg = \sum_{i} \frakg_{i}\) and 
\( \frakh = \sum_{j} \frakh_{j}\) be \(\mZ\)-grade Lie superalgebras. 
Then we have a ``trivial'' superalgebra as follows:
\(\mathfrak{k}_i = \frakg_i\) for \( i< 0\) and 
\(\mathfrak{k}_j = \frakh_j\) for \( j \geqq   0\),  and 
the bracket is defined by 
\[ \SbtE{new}{x}{y} = \begin{cases} \SbtE{\frakg}{x}{y} & \text{if } 
x\in \mathfrak{k}_{i} ,\; 
y\in \mathfrak{k}_{j} (i,j <0)  \\
\SbtE{\frakh}{x}{y} & \text{if  }
x\in \mathfrak{k}_{i} ,\; 
y\in \mathfrak{k}_{j} (i,j \geqq 0)  \\
0 & \text{otherwise.} \end{cases} \] 
%like above proposition. 
\end{prop}

\section{One step extended superalgebra} 
\begin{theorem} Let \(M\) be an \(n\)-dimensional manifold.  
As \eqref{super:cbdl}, we have a superalgebra  
\(\ds \frakh = 
 \Lambda ^{ n } \cbdl{M}\oplus \cdots \oplus  
% \Lambda ^{ 1 } \cbdl{M} \oplus 
 \Lambda ^{ 0 } \cbdl{M} %\oplus \tbdl{M}
\).  
Take \(\ds \ovfrakh = \frakh \oplus \Lambda \tbdl{M}\) with 
\( \ovfrakh_{0} = \tbdl{M}\).  

On \( \ovfrakh \), 
using Lie derivative \(L_{X}\) of \(X \in \tbdl{M}\), we 
define a bracket by  
\begin{align*}
 \SbtE{new}{x}{y} &= \SbtS{x}{y} \text{ for 1-vectors } x,y,   \\    
 \SbtE{new}{x}{y} &= \SbtD{x}{y} \text{ for forms } x,y, \\ 
 \SbtE{new}{x}{y} &= L_{x} y 
 \text{ for 1-vector }  x  \text{  and  form }   y \; \\ 
 \SbtE{new}{y}{x} & = - L_{x} y 
 \text{ for 1-vector }  x  \text{  and  form }   y \; . 
\end{align*}
Then \(\ovfrakh\) becomes a Lie superalgebra.  
\(\frakh\)  is a (super subalgebra) or sub superalgebra of
\(\ds \ovfrakh \) naturally.  

\end{theorem}
\textbf{Proof:} We only check super Jacobi identity by the properties  
\begin{align*}
L_{X} (\alpha \wedge \beta) & = 
L_{X} (\alpha) \wedge \beta  + \alpha \wedge L_{X} (\beta) \\  
L_{\Sbt{X}{Y}}\alpha &= [ L_{X}, L_{Y} ] \alpha =  L_{X}( L_{Y}  \alpha) - 
  L_{Y}( L_{X}  \alpha) \; .  
\end{align*}
In fact, for \(\alpha \in \Lambda^{a} \cbdl{M}\) 
\begin{align*}
\SbtE{new}{\Sbt{X}{Y}} {\alpha} &= 
  L_{X}( L_{Y}  \alpha) - L_{Y}( L_{X}  \alpha) 
= L_{X} \SbtE{new}{Y}{\alpha} - L_{Y} \SbtE{new}{X}{\alpha}
\\&= \SbtE{new}{X}{\SbtE{new}{Y}{\alpha}}
- \SbtE{new}{Y}{\SbtE{new}{X}{\alpha}}
%\\&
= \SbtE{new}{X}{\SbtE{new}{Y}{\alpha}} + \SbtE{new}{\SbtE{new}{X}{\alpha}}{Y}
\\&= \SbtE{new}{X}{\SbtE{new}{Y}{\alpha}}
+ \parity{0 \cdot (a-1)}\SbtE{new}{\SbtE{new}{X}{\alpha}}{Y}
\\
\SbtE{new}{X}{\SbtD{\alpha}{\beta}}
&= L_{X}(\parity{a} d ( \alpha\wedge \beta ))  
%\\& 
= \parity{a} d(L_{X}( \alpha\wedge \beta ))  
= \parity{a} d(L_{X}( \alpha) \wedge \beta 
+  \alpha \wedge L_{X}(\beta) ))  
\\&= \SbtD{ L_{X}\alpha }{\beta} + \SbtD{\alpha}{L_{X}\beta}
= \SbtD{\SbtE{new}{X}{\alpha}}{\beta} + \parity{0\cdot (a-1)} 
\SbtD{\alpha}{\SbtE{new}{X}{\beta}}
\end{align*}
\qed

\subsection{General argument}
Since the degree or weight of \(X\in \tbdl{M}\) is 0, we have 
\begin{equation}
\myOSW{m}{w} = \sum_{k=0}^{n} \mywedge^k \tbdl{M} \mywedge \myCSW{m-k}{w}
\end{equation}
\subsection{About low dimensional Lie algebras}
In this case, 
\begin{align}
\myOSW{m}{w} &= \sum_{k=0}^{n} \mywedge^k \frakg \mywedge \myCSW{m-k}{w}
\\
\sum_{m} \parity{m}
\dim \myOSW{m}{w} &
=\sum_{m}\parity{m} \sum_{k=0}^{n} \tbinom{n}{k}\dim \myCSW{m-k}{w}
=\sum_{m,k}\parity{k}  \tbinom{n}{k} \parity{m-k}\dim \myCSW{m-k}{w}
\\
& =\sum_{k}\parity{k}  \tbinom{n}{k} \sum_{m} \parity{m-k}\dim \myCSW{m-k}{w}
= 0 \notag
\end{align}
Thus, the Euler number is always 0.

\subsection{Double weighted super homology theory on \(\mR^{n}\)}
Let \( x_{1},\ldots, x_{n}\) be a Cartesian coordinate of \(\mR^{n}\). 
Then vector fields and differential forms with 
only polynomial coefficients are written as follows:  
\begin{align*}
X &= \sum_{i=1}^{n} F_{i} (x) \frac{\pdel }{\pdel x_{i}} \quad \text{where }
F_{i}(x) \text{ is a polynomial of } x = ( x_1,\ldots, x_n )
\\
\Omega 
&= \sum_{A} G_{A}(x) dx^{A} \quad \text{where } 
G_{A}(x) \text{ is a polynomial of } x,  
 dx^{A} = dx_{1}^{a_1} \wedge \cdots \wedge dx_{n}^{a_n} \text{ and } 
\end{align*}
\(
A = (a_1,\dots,a_n)\in \{0,1\} ^{n}\), we sometimes denote \(a_{1}+
\cdots+a_{n}\) by \(|A|\).   

They provide Lie sub superalgebras.   
We divide polynomials by homogeneity and define the secondary degree or
weight as follows. 
\begin{definition}
Assume  
\(\ds X = 
 \sum_{i=1}^{n} F_{i} (x) \frac{\pdel }{\pdel x_{i}}\ne 0\). 
The grading in super sense or the weight of \(X\) is 0 by definition. 
If for each \(i\), \( F_{i}(x) =0\) or homogeneous with the same homogeneity
\(h\), i.e., \( F_{i}(tx)= t^{h}F_{i}(x)\) for \(\forall i\),  then we define the secondary weight of \(X\) is \(h-1\). 

Non-zero \(\ds \Omega = \sum_{|A|=a} G_{A}(x) dx^{A}\) has the degree in
super sense or the (primary) weight \( -1 -a\) by the definition. Now the
secondary weight is defined by \(h-1\) when \( G_{A}(tx) = t^{h} G_{A}(x)\)
for \(\forall A\) with \( |A|=a\).  
Let \( \wfg{i,h} \) be the space consisted by 
\((h+1)\)-homogeneous \((-1-i)\)-forms for \( i<0\), and let 
\( \wfg{0,h} \) be the space consisted by 
\((h+1)\)-homogeneous 1-vectors. 
\(\ds \frakg_{i} = \sum_{h} \wfg{i,h}\)  
and 
\(\ds \frakg_{0} = \sum_{h} \wfg{0,h}\) hold.   
\end{definition}

We easily see that 
\begin{prop}
\(\ds \SbtE{new}{ \wfg{i,h} } { \wfg{j,k} } \subset  
 \wfg{i+j,h+k}  
\) for \( i,j \leqq 0\) and \( h,k \geqq -1\).  
\end{prop}
Thus, we have doubly weighted chain complex 
\[\myCSW{m}{w,h}= \sum_{i=1}^{m} \{  \wfg{ w_{i}, h_{i} } \mid \sum_{i} w_{i}
= w \text{ and } \sum_{i} h_{i} = h\}
\;.
\]  
\newcommand{\Eu}{E}
\newcommand{\SW}[1]{SW(#1)\;}
\newcommand{\PW}[1]{PW(#1)\;}
\begin{theorem}
The Betti numbers are 0 when \( w-h \ne 0\) for 
the double weighted homology groups
of the above double weighted chain complex \(\myCSW{m}{w,h}\)
of \(\mR^{n}\).
\end{theorem}
\textbf{Proof:}
Using the Euler vector field 
\(\ds \Eu = \sum_{i=1}^{n} \xb{i} \frac{\pdel}{\pdel \xb{i}} \), we have 
\begin{align*}
\Sbt{\Eu} {x^{P} d x^{A}} &= (|P|+ |A|) x^{P} d x^{A} = 
 (\SW{x^{P} d x^{A}}- \PW{x^{P} d x^{A}} ) x^{P} d x^{A}  \\
\Sbt{\Eu} {x^{P} \pdel/ \pdel  \xb{j}} &
= (|P|- 1) x^{P}\pdel /\pdel\xb{j}
= \SW{x^{P}\pdel /\pdel\xb{j} } x^{P}\pdel /\pdel\xb{j}
\end{align*}
where 
\(
 \PW{\alpha} \) and  
\(\SW{\alpha}\) mean the primary and the secondary weight of \( \alpha \).

Since the double weight of \(\Eu\) is \((0,0)\), 
the definition of the boundary operator implies 
\begin{align*}
\pdel ( \Eu \mw Y_{1} \mw \cdots \mw Y_{m} ) &= - 
\Eu \mw \pdel( Y_{1}  \mw \cdots \mw Y_{m} ) + \sum_{i=1}^{m} 
Y_{1} \mw \cdots \mw \Sbt{\Eu}{Y_{i}}\mw\cdot \mw Y_{m}  
\\\shortintertext{When  \( Y_{i} \in \wfg{w_{i},h_{i}}\), 
\( \Sbt{\Eu}{Y_{i}} = (h_{i}+1) Y_{i} \) 
because \(\Eu\) is the Euler vector field, so we have   
}
\pdel ( \Eu \mw Y_{1} \mw \cdots \mw Y_{m} ) &= - 
\Eu \mw \pdel(Y_{1}\mw \cdots \mw Y_{m} ) + \sum_{i=1}^{m} ( h_{i} - w_{i} ) 
Y_{1} \mw \cdots \mw Y_{m}  
\\ &= - 
\Eu \mw \pdel( Y_{1}  \mw \cdots \mw Y_{m} ) + (h-w)
Y_{1} \mw \cdots \mw Y_{m} 
\end{align*}
for \( Y_{1} \mw \cdots \mw Y_{m} \in \myCSW{m}{w,h}\). 
Thus, if \( h-w \ne 0\), we see each \(m\)-cycle is exact and so the Betti
number is 0. 
\qed

%The discussion of case \(h - w = 0\)  is stacked right now (Nov 08, 2021).  

\kmcomment{
We are interested in what kind of difference 
appear, comparing with \cite{Mik:Miz:super3} where  we have discussed    
the doubly weighted second homology group of the sub superalgebra  
of 
\( \ds 
 \Lambda ^{ 1 } \tbdl{M}\oplus \cdots  \oplus 
 \Lambda ^{ n } \tbdl{M} 
 %\qquad \text{ grade of }\Lambda^{i}\tbdl{M} = i-1 
 \). 
}

\appendix

\newcommand{\ol}[1]{\overline{#1}}
\newcommand{\Sbtn}[2]{\SbtE{n}{#1}{#2}}
\newcommand{\SbtDp}[2]{\SbtDE{ }{#1}{#2}}
\newcommand{\SbtDn}[2]{\SbtDE{+}{#1}{#2}}
\newcommand{\BktR}[2]{\BktD{r}{#1}{#2}}
\newcommand{\BktN}[2]{\BktD{n}{#1}{#2}}
\newcommand{\Sbtr}[2]{\SbtX{r}{#1}{#2}}
% \( \BktR{A}{B} \Sbtr{A}{B} \BktD{u}{A}{B} \)
\section{Appendix}
We have a new super bracket on forms in \cite{Mik:Miz:superForms} which is
older version of this article.  We first fix a candidate as 
\( L(a,b)  d\alpha \we \beta + R(a,b) \alpha \we d\beta \). By super
symmetry and super Jacobi identity we got \(L\) and \(R\) up to
constant. That is our super bracket. 
On the other hand, 
the Schouten bracket is a prototype of super brackets, and there is a way to
understand it.  We first recall the way in the next subsection, then     
in the final subsection, we explain our super bracket can be understood in a similar way. 

\subsection{Schouten bracket}
First we review one procedure to get the Schouten bracket. 
Let \(x,y\) be vector fields, and \(\alpha \) be a differential form on
\(M\). Then we have the magic formula and a formula
\begin{align}
\Lb{x} \alpha &= ( d \iota_{x} + \iota_{x} d ) \alpha \label{magic:formula} \\
\iota_{\Sbt{x}{y}} \alpha &= \Lb{x} \iota_{y} \alpha - \iota_{y} \Lb{x} \alpha\;
,\quad\text{i.e.,}\quad
\iota_{y} \Lb{x} =  \Lb{x} \iota_{y} +  \iota_{\Sbt{y}{x}} 
\label{rare:formula} 
\end{align}

From \eqref{magic:formula}, it is known 
\begin{lemma}
\begin{equation}
 d \circ \inn{x_{1}} \circ \cdots \circ \inn{x_{p}}
- (-1)^{p} \inn{x_{1}}  \circ \cdots \circ \inn{x_{p}} \circ d
 = \sum_{j=1}^{p} (-1)^{j+1}
\inn{x_{1}}\circ \cdots \Lb{x_{j}} \cdots\circ\inn{x_{p}}\; . 
\label{repeated:magic:formula}
\end{equation}
\end{lemma}

Symbolically, we denote them as 
\( \inn{x_{1}} \circ \cdots \circ \inn{x_{p}}  = 
\inn{x_{1}} \we \cdots \we \inn{x_{p}}= X \) and 
\(
 \sum_{j=1}^{p} (-1)^{j+1}
\inn{x_{1}}\circ \cdots \Lb{x_{j}} \cdots\circ\inn{x_{p}}=
 \ol{X}\). 
The relation \eqref{repeated:magic:formula} is symbolically  
\begin{equation} \Sbt{d}{X} =  d \cdot X - \parity{p} X \cdot d = \ol{X}
\end{equation}
We use 
\eqref{rare:formula} a lot and get 
\begin{align*}
\iota_{y}\ol{X} &= \sum_{j=1}^{p} 
\inn{x_{1}}\circ \cdots \circ \inn{y} \Lb{\xb{j}}\circ \cdots\circ\inn{\xb{p}}
= 
 \sum_{j=1}^{p}\inn{x_{1}}\circ \cdots \circ ( \Lb{\xb{j}} \inn{y}  + \inn{ \Sbt{y}{\xb{j}} }
) \circ \cdots\circ\inn{\xb{p}}
\\& =\left( \sum_{j=1}^{p}\parity{p-j}\inn{x_{1}}\circ \cdots \circ (\Lb{\xb{j}}) \circ \cdots\circ 
\inn{\xb{p}} \right) \inn{y} 
+ \sum_{j=1}^{p} \inn{x_{1}}\circ \cdots \circ (  \inn{ \Sbt{y}{\xb{j}} }) \circ \cdots\circ 
\inn{\xb{p}}
\\&= \parity{p-1} \ol{X} \cdot \inn{y} + \Lb{y}{X} \quad\text{where} \quad 
\Lb{y}{X} = \sum_{j=1}^{p} {\xb{1}}\we  \cdots \we  { \Sbt{y}{\xb{j}} }  \we
 \cdots\we \xb{p} \;. 
\\
\inn{z} \iota_{y}\ol{X} &= 
\inn{z} ( \parity{p-1} \ol{X}\cdot \inn{y} + \Lb{y} ) 
=   \parity{p-1} \inn{z}\ol{X}\cdot \inn{y} + \inn{z}\Lb{y}  
=   \parity{2(p-1)} \ol{X}\cdot\inn{z}\cdot \inn{y} 
+ \parity{p-1} \Lb{z} \cdot  \inn{y}
+ \inn{z}\Lb{y} \;.  
\\\shortintertext{
By induction on the degree of \(Y\), we get
}
Y \cdot \ol{X} &= \parity{q(p-1)} \ol{X} \cdot Y + \SbtS{Y}{X}\;,
\quad\text{i.e., } \quad  
\Sbt{Y}{\Sbt{d}{X}} = \SbtS{Y}{X}\;, \\ 
\shortintertext{where}
\SbtS{Y}{X} &= \sum_{k=1}^{q}\parity{(p-1)(q-k)} 
{\yb{1}}\we  \cdots \we  \Lb{\yb{k}} X  \we \cdots\we \yb{q}\\&
= \sum_{j,k} \parity{j+k} \Sbt{\yb{k}}{\xb{j}} \cdot 
\yb{1}\we  \cdots \we  \widehat{\yb{k}}   \we \cdots\we \yb{q} \we  
\xb{1}\we  \cdots \we  \widehat{\xb{j}}   \we \cdots\we \xb{p} 
 \;.
\end{align*}
This \( \SbtS{Y}{X}\) is the Schouten bracket of multi-vector fields \(Y\)
and \(X\).

\subsection{Super bracket on forms}
We apply the explanation so far to reach our bracket on forms. 
Let us 
consider the  \(\mZ\)-graded vector space 
\(\frakh= \sum_{i} \frakhN{i}\) where \(\frakhN{i} =  \cgaiseki{i} \).  
The composition is just the wedge product and it holds   
\[ d \cdot  A = ( d A ) + \parity{a} A \cdot d \quad\text{more precisely  }
\quad 
 d \cdot (A  \cdot \omega) = ( d A )\cdot \omega + \parity{a} A \cdot(d 
\cdot \omega)  \; . 
\]
Unlike the previous discussion,  \(
\Sbt{\Sbt{d}{A}}{B} =  \Sbt{d A} {B} = 0\) because  
\( \Sbt{A}{B} = A \cdot B -  \parity{ab} B \cdot A\) and so  
\( \Sbt{d}{A} = dA\).  
So we try another bracket 
\begin{definition}
\(\Pkt{A}{B} = A \cdot B + \parity{a b} B \cdot A\) and 
\(\BktN{A}{B} = \Pkt{A}{ \Pkt{d}{B}} \). 
\end{definition}
Direct computation implies 
\begin{align}
\BktN{A}{B} &= \Pkt{A}{ \Pkt{d}{B}} 
= 2 A \cdot (d B) + 2 \parity{a+b+ab} B \cdot (dA) + 4 \parity{b} A \cdot B \cdot d 
\label{extra:Sbtn}
\end{align}
Since \(d\cdot d=0\), \(d \cdot \beta = d\beta + \parity{b} \beta \cdot
d\) and \( d\beta \in \frakh\), the algebra generated \(d\) and \(\frakh\)
is \( \frakh \oplus \frakh \cdot d\). 
Applying the projection from 
\( \frakh \oplus \frakh \cdot d\) 
onto \(\frakh\) to \eqref{extra:Sbtn}, we get and  
define a revised bracket as follows.
\begin{definition}
\begin{equation} \label{eqn:sbtr:eg}
\BktR{A}{B} = 
  A \cdot (d B) +  \parity{a+b+ab} B \cdot (dA) = \parity{a} d( A \cdot B ) 
\end{equation}
for each \(A \in \frakhN{a}\) and \( B \in \frakhN{b}\). 
\end{definition}
\begin{theorem}
The bracket \(\BktR{\cdot}{\cdot}\) is a superalgebra
bracket by the new grading \(\ndeg{\frakhN{i}} = i+1\). 
\end{theorem}

\nocite{Mik:Miz:homogPoisson} 
%\cite{Mik:Miz:super1} 
\nocite{Mik:Miz:super2} 
\nocite{Mik:Miz:super3} 
\nocite{Mik:Miz:superLowDim}

\bibliographystyle{plain}
\bibliography{km_refs}
%\printbibliography

\end{document}